%
%
%
%


\magnification=1200
\pretolerance=500 \tolerance=1000 \brokenpenalty=5000
\hsize=12.5cm   
\vsize=19cm
\hoffset=0.4cm
\voffset=1cm   
\parskip3pt plus 1pt
\parindent=0.6cm
\let\sl=\it
\def\\{\hfil\break}


\font\seventeenbf=cmbx10 at 17.28pt

\font\twelvebf=cmbx10 at 12pt
\font\eightbf=cmbx8
\font\sevenbf=cmbx10 at 7pt
\font\sixbf=cmbx6

\font\eighti=cmmi8
\font\sixi=cmmi6

\font\eightrm=cmr8
\font\sixrm=cmr6

\font\eightsy=cmsy8
\font\sixsy=cmsy6

\font\eightit=cmti8
\font\eighttt=cmtt8
\font\eightsl=cmsl8

\font\seventeenbsy=cmbsy10 at 17.28pt

\font\twelvebsy=cmbsy10 at 12pt
\font\tenbsy=cmbsy10
\font\eightbsy=cmbsy8
\font\sevenbsy=cmbsy7
\font\sixbsy=cmbsy6
\font\fivebsy=cmbsy5

\font\tenmsa=msam10

\font\sevenmsa=msam7
\font\fivemsa=msam5
\newfam\msafam
  \textfont\msafam=\tenmsa
  \scriptfont\msafam=\sevenmsa
  \scriptscriptfont\msafam=\fivemsa

\font\tenmsb=msbm10
\font\eightmsb=msbm8
\font\sevenmsb=msbm7
\font\fivemsb=msbm5
\newfam\msbfam
  \textfont\msbfam=\tenmsb
  \scriptfont\msbfam=\sevenmsb
  \scriptscriptfont\msbfam=\fivemsb
\def\Bbb{\fam\msbfam\tenmsb}

\font\tenCal=eusm10
\font\sevenCal=eusm7
\font\fiveCal=eusm5
\newfam\Calfam
  \textfont\Calfam=\tenCal
  \scriptfont\Calfam=\sevenCal
  \scriptscriptfont\Calfam=\fiveCal
\def\Cal{\fam\Calfam\tenCal}

\font\teneuf=eusm10
\font\teneuf=eufm10
\font\seveneuf=eufm7
\font\fiveeuf=eufm5
\newfam\euffam
  \textfont\euffam=\teneuf
  \scriptfont\euffam=\seveneuf
  \scriptscriptfont\euffam=\fiveeuf
\def\euf{\fam\euffam\teneuf}

\font\seventeenbfit=cmmib10 at 17.28pt

\font\twelvebfit=cmmib10 at 12pt
\font\tenbfit=cmmib10
\font\eightbfit=cmmib8
\font\sevenbfit=cmmib7
\font\sixbfit=cmmib6
\font\fivebfit=cmmib5
\newfam\bfitfam
  \textfont\bfitfam=\tenbfit
  \scriptfont\bfitfam=\sevenbfit
  \scriptscriptfont\bfitfam=\fivebfit


\catcode`\@=11
\def\eightpoint{%
  \textfont0=\eightrm \scriptfont0=\sixrm \scriptscriptfont0=\fiverm
  \def\rm{\fam\z@\eightrm}%
  \textfont1=\eighti \scriptfont1=\sixi \scriptscriptfont1=\fivei
  \def\oldstyle{\fam\@ne\eighti}%
  \textfont2=\eightsy \scriptfont2=\sixsy \scriptscriptfont2=\fivesy
  \textfont\itfam=\eightit
  \def\it{\fam\itfam\eightit}%
  \textfont\slfam=\eightsl
  \def\sl{\fam\slfam\eightsl}%
  \textfont\bffam=\eightbf \scriptfont\bffam=\sixbf
  \scriptscriptfont\bffam=\fivebf
  \def\bf{\fam\bffam\eightbf}%
  \textfont\ttfam=\eighttt
  \def\tt{\fam\ttfam\eighttt}%
  \textfont\msbfam=\eightmsb
  \def\Bbb{\fam\msbfam\eightmsb}%
  \abovedisplayskip=9pt plus 2pt minus 6pt
  \abovedisplayshortskip=0pt plus 2pt
  \belowdisplayskip=9pt plus 2pt minus 6pt
  \belowdisplayshortskip=5pt plus 2pt minus 3pt
  \smallskipamount=2pt plus 1pt minus 1pt
  \medskipamount=4pt plus 2pt minus 1pt
  \bigskipamount=9pt plus 3pt minus 3pt
  \normalbaselineskip=9pt
  \setbox\strutbox=\hbox{\vrule height7pt depth2pt width0pt}%
  \let\bigf@ntpc=\eightrm \let\smallf@ntpc=\sixrm
  \normalbaselines\rm}
\catcode`\@=12

\def\eightpointbf{%
 \textfont0=\eightbf   \scriptfont0=\sixbf   \scriptscriptfont0=\fivebf
 \textfont1=\eightbfit \scriptfont1=\sixbfit \scriptscriptfont1=\fivebfit
 \textfont2=\eightbsy  \scriptfont2=\sixbsy  \scriptscriptfont2=\fivebsy
 \eightbf
 \baselineskip=10pt}

\def\tenpointbf{%
 \textfont0=\tenbf   \scriptfont0=\sevenbf   \scriptscriptfont0=\fivebf
 \textfont1=\tenbfit \scriptfont1=\sevenbfit \scriptscriptfont1=\fivebfit
 \textfont2=\tenbsy  \scriptfont2=\sevenbsy  \scriptscriptfont2=\fivebsy
 \tenbf}
        
\def\twelvepointbf{%
 \textfont0=\twelvebf   \scriptfont0=\eightbf   \scriptscriptfont0=\sixbf
 \textfont1=\twelvebfit \scriptfont1=\eightbfit \scriptscriptfont1=\sixbfit
 \textfont2=\twelvebsy  \scriptfont2=\eightbsy  \scriptscriptfont2=\sixbsy
 \twelvebf
 \baselineskip=14.4pt}

\def\seventeenpointbf{%
 \textfont0=\seventeenbf  \scriptfont0=\twelvebf  \scriptscriptfont0=\eightbf
 \textfont1=\seventeenbfit\scriptfont1=\twelvebfit\scriptscriptfont1=\eightbfit
 \textfont2=\seventeenbsy \scriptfont2=\twelvebsy \scriptscriptfont2=\eightbsy
 \seventeenbf
 \baselineskip=20.736pt}
 

\newdimen\srdim \srdim=\hsize
\newdimen\irdim \irdim=\hsize
\def\NOSECTREF#1{\noindent\hbox to \srdim{\null\dotfill ???(#1)}}
\def\SECTREF#1{\noindent\hbox to \srdim{\csname REF\romannumeral#1\endcsname}}
\def\INDREF#1{\noindent\hbox to \irdim{\csname IND\romannumeral#1\endcsname}}
\newlinechar=`\^^J
\def\openauxfile{
  \immediate\openin1\jobname.
  \ifeof1
  \let\sectref=\NOSECTREF \let\indref=\NOSECTREF
  \else
  \input \jobname.aux
  \let\sectref=\SECTREF \let\indref=\INDREF
  \fi
  \immediate\openout1=\jobname.aux
  \let\END=\end \def\end{\immediate\closeout1\END}}
        
\newbox\titlebox   \setbox\titlebox\hbox{\hfil}
\newbox\sectionbox \setbox\sectionbox\hbox{\hfil}
\def\folio{\ifnum\pageno=1 \hfil \else \ifodd\pageno
           \hfil {\eightpoint\copy\sectionbox\kern8mm\number\pageno}\else
           {\eightpoint\number\pageno\kern8mm\copy\titlebox}\hfil \fi\fi}
\footline={\hfil}
\headline={\folio}

\def\titlerunning#1{\setbox\titlebox\hbox{\eightpoint #1}}
\def\title#1{\noindent\hfil$\smash{\hbox{\seventeenpointbf #1}}$\hfil
             \titlerunning{#1}\medskip}

\newcount\numbersection \numbersection=-1
\def\sectionrunning#1{\setbox\sectionbox\hbox{\eightpoint #1}
  \immediate\write1{\string\def \string\REF 
      \romannumeral\numbersection \string{%
      \noexpand#1 \string\dotfill \space \number\pageno \string}}}
\def\section#1{%
  \par\vskip0.666cm\penalty -100
  \vbox{\baselineskip=14.4pt\noindent{{\twelvepointbf #1}}}
  \vskip2pt
  \penalty 500
  \advance\numbersection by 1
  \sectionrunning{#1}}

\def\subsection#1|{%
  \par\vskip0.5cm\penalty -100
  \vbox{\noindent{{\tenpointbf #1}}}
  \vskip1pt
  \penalty 500}

\newcount\numberindex \numberindex=0  
\def\index#1#2{%
  \advance\numberindex by 1
  \immediate\write1{\string\def \string\IND #1%
     \romannumeral\numberindex \string{%
     \noexpand#2 \string\dotfill \space \string\S \number\numbersection, 
     p.\string\ \space\number\pageno \string}}}

\newdimen\itemindent \itemindent=\parindent

\def\item#1{\par\noindent\hangindent\itemindent%
            \rlap{#1}\kern\itemindent\ignorespaces}
\def\itemitem#1{\par\noindent\hangindent2\itemindent%
            \kern\itemindent\rlap{#1}\kern\itemindent\ignorespaces}
\def\itemitemitem#1{\par\noindent\hangindent3\itemindent%
            \kern2\itemindent\rlap{#1}\kern\itemindent\ignorespaces}

\long\def\claim#1|#2\endclaim{\par\vskip 5pt\noindent 
{\tenpointbf #1.}\ {\sl #2}\par\vskip 5pt}

\def\proof{\noindent{\sl Proof}}
\def\dummy{\phantom{$\hat{\hat{\hat{\hat f}}}$}}

\def\today{\ifcase\month\or
January\or February\or March\or April\or May\or June\or July\or August\or
September\or October\or November\or December\fi \space\number\day,
\number\year}

\catcode`\@=11
\newcount\@tempcnta \newcount\@tempcntb 
\def\timeofday{{%
\@tempcnta=\time \divide\@tempcnta by 60 \@tempcntb=\@tempcnta
\multiply\@tempcntb by -60 \advance\@tempcntb by \time
\ifnum\@tempcntb > 9 \number\@tempcnta:\number\@tempcntb
  \else\number\@tempcnta:0\number\@tempcntb\fi}}
\catcode`\@=12

\def\bibitem#1&#2&#3&#4&%
{\hangindent=1.66cm\hangafter=1
\noindent\rlap{\hbox{\eightpointbf #1}}\kern1.66cm{\rm #2}{\sl #3}{\rm #4.}} 


\def\bC{{\Bbb C}}
\def\bG{{\Bbb G}}

\def\bN{{\Bbb N}}
\def\bP{{\Bbb P}}
\def\bQ{{\Bbb Q}}

\def\bZ{{\Bbb Z}}
\def\bone{{\mathchoice {\rm 1\mskip-4.2mu l} {\rm 1\mskip-4mu l}
{\rm 1\mskip-4.5mu l} {\rm 1\mskip-5mu l}}}

\def\gm{{\euf m}}


\def\cK{{\Cal K}}

\def\cO{{\Cal O}}

\def\cI{{\Cal I}}
\def\cJ{{\Cal J}}


\def\bu{{\scriptstyle\bullet}}

\def\square{{\hfill \hbox{
\vrule height 1.453ex  width 0.093ex  depth 0ex
\vrule height 1.5ex  width 1.3ex  depth -1.407ex\kern-0.1ex
\vrule height 1.453ex  width 0.093ex  depth 0ex\kern-1.35ex
\vrule height 0.093ex  width 1.3ex  depth 0ex}}}
\def\qed{\phantom{$\quad$}\hfill$\square$\medskip}
\def\hexnbr#1{\ifnum#1<10 \number#1\else
 \ifnum#1=10 A\else\ifnum#1=11 B\else\ifnum#1=12 C\else
 \ifnum#1=13 D\else\ifnum#1=14 E\else\ifnum#1=15 F\fi\fi\fi\fi\fi\fi\fi}
\def\msatype{\hexnbr\msafam}
\def\msbtype{\hexnbr\msbfam}
\mathchardef\restriction="3\msatype16   
\mathchardef\compact="3\msatype62
\mathchardef\smallsetminus="2\msbtype72   \let\ssm\smallsetminus
\mathchardef\subsetneq="3\msbtype28
\mathchardef\supsetneq="3\msbtype29
\mathchardef\leqslant="3\msatype36   \let\le\leqslant
\mathchardef\geqslant="3\msatype3E   \let\ge\geqslant
\mathchardef\ltimes="2\msbtype6E
\mathchardef\rtimes="2\msbtype6F

\let\ol=\overline

\let\wt=\widetilde

\let\text=\hbox
\def\build#1^#2_#3{\mathrel{\mathop{\null#1}\limits^{#2}_{#3}}}
\def\mertorelbar{\vrule width0.6ex height0.65ex depth-0.55ex}
\def\merto{\mathrel{\mertorelbar\kern1.3pt\mertorelbar\kern1.3pt\mertorelbar
    \kern1.3pt\mertorelbar\kern-1ex\raise0.28ex\hbox{${\scriptscriptstyle>}$}}}
\let\lra=\longrightarrow
\catcode`\@=11
\newdimen\@rrowlength \@rrowlength=6ex
\def\ssrelbar{\vrule width\@rrowlength height0.64ex depth-0.56ex\kern-4pt}
\def\llra#1{\@rrowlength=#1\ssrelbar\rightarrow}
\catcode`\@=12

\def\card{\mathop{\rm card}\nolimits}
\def\Div{\mathop{\rm Div}\nolimits}

\def\Tr{\mathop{\rm Tr}\nolimits}

\def\Hom{\mathop{\rm Hom}\nolimits}

\def\Pic{\mathop{\rm Pic}\nolimits}

\def\Sing{\mathop{\rm Sing}\nolimits}

\def\rank{\mathop{\rm rank}\nolimits}
\def\lcm{\mathop{\rm lcm}\nolimits}
\def\deg{\mathop{\rm deg}\nolimits}

\def\Vol{\mathop{\rm Vol}\nolimits}

\def\dbar{{\overline\partial}}
\def\ddbar{{\partial\overline\partial}}


\def\GG{{\rm GG}}

\long\def\InsertFig#1 #2 #3 #4\EndFig{\par
\hbox{\hskip #1mm$\vbox to#2mm{\vfil\special{" 
(/home/demailly/psinputs/grlib.ps) run
#3}}#4$}}
\long\def\LabelTeX#1 #2 #3\ELTX{\rlap{\kern#1mm\raise#2mm\hbox{#3}}}


\openauxfile

\title{Holomorphic Morse inequalities and}
\smallskip
\title{the Green-Griffiths-Lang conjecture}
\titlerunning{Holomorphic Morse inequalities and
the Green-Griffiths-Lang conjecture}
\medskip
\centerline{\twelvebf Jean-Pierre Demailly}
\medskip
\centerline{Universit\'e de Grenoble I, D\'epartement de Math\'ematiques}
\centerline{Institut Fourier, 38402 Saint-Martin d'H\`eres, France}
\centerline{{\it e-mail\/}: {\tt demailly@fourier.ujf-grenoble.fr}}
\vskip15pt
\line{\hfill\it Dedicated to the memory of Eckart Viehweg}
\vskip15pt

\noindent
{\bf Abstract.} 
The goal of this work is to study the existence and properties of non
entire curves $f:\bC\to X$ drawn in a complex irreducible
$n$-dimensional variety~$X$, and more specifically to show that they 
must satisfy certain global algebraic or differential equations 
as soon as $X$ is projective of general type. By means of
holomorphic Morse inequalities and a probabilistic analysis of
the cohomology of jet spaces, we are able to prove a significant step 
of a generalized version of the Green-Griffiths-Lang conjecture on the
algebraic degeneracy of entire curves.
\smallskip

\noindent
{\bf R\'esum\'e.} Le but de ce travail est d'\'etudier l'existence et les 
propri\'et\'es des courbes enti\`eres $f:\bC \to X$ trac\'ees sur une 
var\'et\'e complexe 
irr\'eductible de dimension~$n$, et plus pr\'ecis\'ement de montrer que ces
courbes doivent satisfaire \`a certaines \'equations alg\'ebriques 
ou diff\'erentielles globales d\`es que $X$ est projective
de type g\'en\'eral. Au moyen des in\'egalit\'es de Morse holomorphes et 
d'une analyse probabiliste de la cohomologie des espaces de jets, nous 
d\'emontrons une premi\`ere \'etape significative d'une version 
g\'en\'eralis\'ee de la conjecture de Green-Griffiths-Lang sur la 
d\'eg\'enerescence alg\'ebrique des courbes enti\`eres.
\bigskip

\noindent
{\bf Key words.}
Chern curvature, holomorphic Morse inequality, jet bundle, cohomology group, 
entire curve, algebraic degeneration, weighted projective space,
Green-Griffiths-Lang conjecture
\smallskip
\noindent
{\bf Mots-cl\'es.}
Courbure de Chern, in\'egalit\'e de Morse holomorphe, fibr\'e de jets, 
groupe de cohomologie,  courbe enti\`ere, d\'eg\'en\'erescence alg\'ebrique, 
espace projectif \`a poids, conjecture de Green-Griffiths-Lang.
\medskip

\noindent
{\bf MSC 2010 Classification.} 32Q45, 32L20, 14C30

\section{0. Introduction}

Let $X$ be a complex $n$-dimensional manifold$\,$; most of the time we will 
assume that $X$ is compact and even projective algebraic. By an ``entire
curve'' we always mean a non constant holomorphic map defined on the
whole complex line~$\bC$, and we say that it is algebraically degenerate
if its image is contained in a proper algebraic subvariety of the ambient
variety. If~$\mu:\smash{\wt X}\to X$ is a modification and $f:\bC\to X$ is
an entire curve whose image
$f(\bC)$ is not contained in the image $\mu(E)$ of the exceptional
locus, then $f$ admits a unique lifting $\smash{\wt f}:\bC\to\smash{\wt X}$.
For this reason, the study of the algebraic degeneration of $f$ is
a birationally invariant problem, and singularities do not play an essential 
role at this stage. We~will therefore assume that $X$ is non singular, 
possibly after performing a suitable composition of blow-ups. 
We are interested more generally in the 
situation where the tangent bundle $T_X$ is equipped with a  
{\it linear subspace} $V\subset T_X$, that is, an irreducible 
complex analytic subset of the total space of $T_X$ such that
\smallskip\noindent
(0.1)~~all fibers $V_x:=V\cap T_{X,x}$ are vector subspaces of $T_{X,x}$.
\smallskip\noindent
Then the problem is to study entire curves $f:\bC\to X$ which are tangent
to $V$, i.e.\ such that $f_*T_\bC\subset V$. We will refer to a pair 
$(X,V)$ as being a {\it directed variety} (or
{\it directed manifold}). A morphism of directed varieties
$\Phi:(X,V)\to(Y,W)$ is a holomorphic map $\Phi:X\to Y$ such that
$\Phi_*V\subset W\,$; by the irreducibility, it is enough to check this
condition over the dense open subset $X\ssm \Sing(V)$
where $V$ is actually a subbundle. Here $\Sing(V)$ denotes the indeterminacy set
of the associated meromorphic map $\alpha:X\merto G_r(T_X)$ to the 
Grassmannian bbundle of $r$-planes in $T_X$, $r=\rank V\,$; we thus have 
$V_{|X\ssm \Sing(V)}=\alpha^*S$ where $S\to G_r(T_X)$ is the tautological 
subbundle of $G_r(T_X)$. In that way, we get a category, 
and we will be mostly interested in the subcategory whose objects $(X,V)$ 
are projective algebraic manifolds equipped with algebraic linear
subspaces. Notice that an entire curve $f:\bC\to X$ tangent to $V$
is just a morphism $f:(\bC,T_\bC)\to (X,V)$.

The case where $V=T_{X/S}$ is the relative tangent space of some
fibration $X\to S$ is of special interest, and so is the
case of a foliated variety (this is the situation where the sheaf of
sections $\cO(V)$ satisfies the Frobenius integrability condition
$[\cO(V),\cO(V)]\subset\cO(V)$); however, it is very useful to 
allow as well non integrable linear subspaces $V$. We refer to 
$V=T_X$ as being the {\it absolute case}. Our main target
is the following deep conjecture concerning the algebraic degeneracy 
of entire curves, which generalizes similar statements made in
[GG79] (see also [Lang86, Lang87]).

\claim (0.2) Generalized Green-Griffiths-Lang conjecture|Let $(X,V)$ be a 
projective directed manifold such that the canonical sheaf $K_V$ is big 
$($in the absolute case $V=T_X$, this means that $X$ is a variety of
general type, and in the relative case we will say that $(X,V)$ is of
general type$)$. Then there should exist an algebraic subvariety 
$Y\subsetneq X$ such that every non constant entire curve $f:\bC\to X$ 
tangent to $V$ is contained in~$Y$.
\endclaim

The precise meaning of $K_V$ and of its bigness will be explained
below \hbox{--}~our definition {\it does not coincide} with other
frequently used definitions and is in our view better suited to
the study of entire curves of $(X,V)$.  One says that $(X,V)$ is
Brody-hyperbolic when there are no entire curves tangent to $V$.
According to (generalized versions of) conjectures of Kobayashi
[Kob70, Kob76] the hyperbolicity of $(X,V)$ should imply that $K_V$ is
big, and even possibly ample, in a suitable sense. It~would then
follow from conjecture (0.2) that $(X,V)$ is hyperbolic if and only if
for every irreducible variety $Y\subset X$, the linear subspace
$V_{\wt Y}=\overline{\smash{T_{\wt Y\ssm E}}
  \cap\mu_*\kern-4pt{}^{-1}V}\subset T_{\wt Y}$ has a big canonical
sheaf whenever $\mu:\smash{\wt Y}\to Y$ is a desingularization and $E$
is the exceptional locus.

The most striking fact known at this date on the Green-Griffiths-Lang 
conjecture is a recent result of Diverio, Merker and Rousseau [DMR10]
in the absolute case, confirming the statement when 
$X\subset\bP^{n+1}_\bC$ is a generic non
singular hypersurface of large degree~$d$, with a (non optimal) sufficient
lower bound $d\ge\smash{2^{n^5}}$.
Their proof is based in an essential way on a strategy developed by Siu 
[Siu02, Siu04], combined with techniques of [Dem95]. Notice that if the 
Green-Griffiths-Lang conjecture holds true, a much stronger and probably
optimal result would be true, namely all smooth hypersurfaces of degree 
$d\ge n+3$
would satisfy the expected algebraic degeneracy statement. Moreover, by results
of Clemens [Cle86] and Voisin [Voi96], a (very) generic hypersurface of degree 
$d\ge 2n+1$ would in fact be hyperbolic for every $n\ge 2$. Such a generic
hyperbolicity statement has been obtained unconditionally 
by McQuillan [McQ98, McQ99] when $n=2$ and $d\ge 35$, and by 
Demailly-El Goul [DEG00] when $n=2$ and $d\ge 21$. 
Recently Diverio-Trapani [DT10] proved the same result when $n=3$ and 
$d\ge 593$. By defi\-nition, proving the algebraic degeneracy means finding
a non zero polynomial $P$ on $X$ such that all entire curves $f:\bC\to X$
satisfy $P(f)=0$. All known methods of proof are based on establishing 
first the existence of certain algebraic differential equations 
$P(f\,;\,f',f'',\ldots,f^{(k)})=0$ of some order $k$, and then trying to find
enough such equations so that they cut out a proper algebraic
locus $Y\subsetneq X$.

Let $J_kV$ be the space of $k$-jets of curves $f:(\bC,0)\to X$ tangent
to $V$. One defines the sheaf $\cO(E^\GG_{k,m}V^*)$ of jet differentials of
order $k$ and degree $m$ to be the sheaf of holomorphic functions 
$P(z;\xi_1,\ldots\xi_k)$ on $J_kV$ which are homogeneous polynomials 
of degree $m$ on the fibers of $J^kV\to X$ with respect to local
coordinate derivatives $\xi_j=f^{(j)}(0)$ (see below in case $V$ has
singularities). The degree $m$ considered here is the weighted 
degree with respect to the natural $\bC^*$ action on $J^kV$ defined
by $\lambda\cdot f(t):=f(\lambda t)$, i.e.\ by reparametrizing the curve
with a homothetic change of variable. Since $(\lambda\cdot f)^{(j)}(t)=
\lambda^jf^{(j)}(\lambda t)$, the weighted action is given in coordinates
by
$$
\lambda\cdot(\xi_1,\xi_2,\ldots,\xi_k)=
(\lambda\xi_1,\lambda^2\xi_2,\ldots,\lambda^k\xi_k).\leqno(0.3)
$$
One of the major tool of the theory is the following result due
to Green-Griffiths [GG79] (see also [Blo26], 
[Dem95, Dem97], [SY96a, SY96b], [Siu97]).

\claim (0.4) Fundamental vanishing theorem|Let $(X,V)$ be a directed projective
variety and $f:(\bC,T_\bC)\to (X,V)$ an entire curve tangent to $V$.
Then for every global section $P\in H^0(X,\smash{E^\GG_{k,m}V^*}
\otimes \cO(-A))$ where $A$ is an ample divisor of $X$, one has 
$P(f\,;\,f',f'',\ldots,f^{(k)})=0$.
\endclaim

It is expected that the global sections of 
$H^0(X,E^\GG_{k,m}V^*\otimes \cO(-A))$ are precisely those which ultimately
define the algebraic locus $Y\subsetneq X$ where the curve $f$ should lie.
The problem is then reduced to the question of showing that there are many
non  zero sections of $H^0(X,E^\GG_{k,m}V^*\otimes \cO(-A))$, and further,
understanding what is their joint base locus. The first part of this program
is the main result of the present paper.

\claim (0.5) Theorem|Let $(X,V)$ be a directed projective variety
such that $K_V$ is big and let $A$ be an ample divisor. Then 
for $k\gg 1$ and $\delta\in\bQ_+$ small enough, $\delta\le c(\log k)/k$,
the number of sections
$h^0(X,\smash{E^\GG_{k,m}V^*}\otimes\cO(-m\delta A))$ has maximal growth, 
i.e.\ is larger that $c_km^{n+kr-1}$ for some~$m\ge m_k$, 
where $c,\,c_k>0$, $n=\dim X$ and $r=\rank V$. In particular, entire curves
$f:(\bC,T_\bC)\to(X,V)$ satisfy $($many$)$ algebraic differential 
equations.
\endclaim

The statement is very elementary to check when $r=\rank V=1$, and 
therefore when $n=\dim X=1$. In higher dimensions  $n\ge 2$, only 
very partial results were known at this point, 
concerning merely the absolute case $V=T_X$. In dimension~$2$,
Theorem~0.5 is a consequence of the Riemann-Roch calculation of 
Green-Griffiths [GG79], 
combined with a vanishing theorem due to Bogomolov [Bog79] -- the latter 
actually only applies to the top cohomology group $H^n$, and things 
become much more delicate when extimates of intermediate cohomology
groups are needed. In higher dimensions, Diverio [Div09] proved the 
existence of sections of $H^0(X,\smash{E^\GG_{k,m}V^*}\otimes\cO(-1))$
whenever $X$ is a hypersurface of $\bP^{n+1}_\bC$ of high degree $d\ge
d_n$, assuming $k\ge n$ and $m\ge m_n$. More recently, Merker [Mer10]
was able to treat the case of arbitrary hypersurfaces of general type,
i.e.\ $d\ge n+3$, assuming this time $k$ to be very large.  The latter
result is obtained through explicit algebraic calculations of the
spaces of sections, and the proof is computationally very
intensive. B\'erczi [Ber10] also obtained related results with a
different approach based on residue formulas, assuming $d\ge 2^{7n\log
  n}$.

All these approaches are algebraic in nature, and while they use some form of
holomorphic Morse inequalities [Dem85], they only require a very special
elementary algebraic case, namely the lower bound 
$$
h^0(X,L^{\otimes m})\ge {m^n\over n!}(A^n-n\,A^{n-1}\cdot B)-o(m^n)
$$
for $L=\cO(A-B)$ with $A,\,B$ nef (cf.\ Trapani [Tra95]). Here,
our techniques are based on more elaborate curvature estimates in the spirit
of Cowen-Griffiths [CG76]. They require the stronger analytic form 
of holomorphic Morse inequalities.

\claim (0.6) Holomorphic Morse inequalities {\rm ([Dem85])}|Let $(L,h)$
be a holomorphic line bundle on a compact complex manifold $X$, equipped
with a smooth hermitian metric $h$, and let $E$ be a holomorphic 
vector bundle. Denote by $\Theta_{L,h}=-{i\over 2\pi}\ddbar\log h$ the curvature
form of $(L,h)$ and consider the open set
$$
X(L,h,q)=\big\{x\in X\,;\;\Theta_{L,h}(x)~\hbox{has signature $(n-q,q)$}\big\}
\quad\hbox{$(q$-index set of $\Theta_{L,h})$},
$$
so that there is a partition 
$X=S\cup\bigcup_{0\le q\le n}\!X(L,h,q)$ where 
\hbox{$S=\{\det\Theta_{L,h}(x){=}0\}$}.
Then, if we put $r=\rank E$, we have asymptotically
as $m$ tends to infinity$\,:$
\smallskip
\item{\rm (a)} $($Weak Morse inequalities$)$
$$
h^q(X, E\otimes L^{\otimes m})\le {m^n\over n!}\,r\int_{X(L,h,q)}
(-1)^q\Theta_{L,h}^n+o(m^n).
$$
\item{\rm (b)} $($Strong Morse inequalities$)$ If
$X(L,h,\le q)=\coprod_{j\le q}X(L,h,j)$, then
$$
\sum_{j=0}^q(-1)^{q-j}h^j(X, E\otimes L^{\otimes m})\le 
{m^n\over n!}\,r\int_{X(L,h,\le q)}(-1)^q\Theta_{L,h}^n+o(m^n).
$$
\item{\rm (c)} $($Lower bound on $h^0)$
$$
h^0(X, E\otimes L^{\otimes m})-h^1(X, E\otimes L^{\otimes m})\ge 
{m^n\over n!}\,r\int_{X(L,h,\le 1)}\Theta_{L,h}^n-o(m^n).
$$
\vskip-5pt
\endclaim

The proof of the above is based on refined spectral estimates for the complex
Laplace-Beltrami operators. Observe that (0.6~c) is just the special
case of (0.6~b) when $q=1$. It has been recently observed that these
inequalities should be optimal in the sense that the asymptotic
cohomology functional $\smash{\widehat h}^q(X,L):=\limsup_{m\to+\infty}{n!\over
m^n}h^0(X,L^{\otimes m})$ satisfies
$$
\widehat h^q(X,L)\le \inf_{h\in C^\infty}\int_{X(L,h,q)}(-1)^q\Theta_{L,h}^n,
\leqno(0.7)
$$
and that conjecturally the inequality should be an equality; it is proved in
[Dem10a], [Dem10b] that this is indeed the case if $n\le 2$ or $q=0$, at least
when $X$ is projective algebraic.

Notice that holomorphic Morse inequalities are essentially
insensitive to singularities, as we can pass to non singular models
and blow-up $X$ as much as we want: if $\mu:\smash{\wt X}\to X$ is a
modification then $\mu_*\cO_{\wt X}=\cO_X$ and $R^q\mu_*\cO_{\wt X}$
is supported on a codimension $1$ analytic subset (even codimension $2$ if
$X$ is smooth).  It follows by the
Leray spectral sequence that the estimates for $L$ on $X$ or for
$\smash{\wt L}=\mu^*L$ on $\smash{\wt X}$ differ by negligible terms
$O(m^{n-1})$. Finally, we can even work with singular hermitian metrics $h$
which have analytic singularities with positive rational coefficients, 
that is, one can write
locally $h=e^{-\varphi}$ where, possibly after blowing up, 
$$
\varphi(z)=c\log\sum_j|g_j|^2~~\hbox{mod $C^\infty$, with $c\in\bQ_+$ and
$g_j$ holomorphic}.\leqno(0.8)
$$
Especially, $\varphi$ is smooth on some Zariski open set $X\ssm Z$ where
$Z=\bigcap g_j^{-1}(0)$, and it has logarithmic poles along $Z$.
Blowing-up the ideal sheaf $\cJ=(g_j)$ leads to divisorial
singularities, and then by replacing $L$ with $\smash{\wt L}
=\mu^*L\otimes\cO(-E)$ 
where $E\in\Div_\bQ(\smash{\wt X})$ is the singularity divisor, we see that
holomorphic Morse inequalities still hold for the sequence of groups
$H^q(X,E\otimes L^{\otimes m}\otimes\cI(h^{\otimes m}))$ where
$\cI(h^{\otimes m})$ is the multiplier ideal sheaf of~$h^{\otimes m}$
(see Bonavero [Bon93] for more details). In the case of linear
subspaces $V\subset T_X$, we introduce singular hermitian metrics as
follows.

\claim (0.9) Definition|A singular hermitian metric on a linear
subspace $V\subset T_X$ is a metric $h$ on the fibers of $V$ such that
the function $\log h:\xi\mapsto\log|\xi|_h^2$ is locally integrable 
on the total space of $V$.
\endclaim

Such a metric can also be viewed as a singular hermitian metric on
the tauto\-logical line bundle $\cO_{P(V)}(-1)$ on the projectivized bundle
$P(V)=V\ssm\{0\}/\bC^*$, and therefore its dual metric $h^*$ defines a 
curvature current $\Theta_{\cO_{P(V)}(1),h^*}$ of type $(1,1)$ on
$P(V)\subset P(T_X)$, such that
$$p^*\Theta_{\cO_{P(V)}(1),h^*} ={i\over 2\pi}\ddbar\log h,\qquad
\hbox{where $p:V\ssm\{0\}\to P(V)$}.
$$
If $\log h$ is quasi-plurisubharmonic $($or quasi-psh, which means psh 
modulo addition of a smooth function$)$ on $V$, then $\log h$ is indeed
locally integrable, and we have moreover
$$\Theta_{\cO_{P(V)}(1),h^*}\ge -C\omega\leqno(0.10)$$
for some smooth positive $(1,1)$-form on $P(V)$ and some constant $C>0\;$;
conversely, if (0.10) holds, then $\log h$ is quasi psh.

\claim (0.11) Definition|We will say that a singular hermitian metric $h$ on
$V$ is {\rm admissible} if $h$ can be written as $h=e^\varphi h_{0|V}$ where
$h_0$ is a smooth positive definite hermitian on $T_X$ and $\varphi$
is a quasi-psh weight with analytic singularities on $X$, as in $(0.9)$. 
Then $h$ can be seen as a singular hermitian metric on $\cO_{P(V)}(1)$,
with the property that it induces a smooth positive definite metric
on a Zariski open set $X'\subset X\ssm \Sing(V)\,;$ we will denote by
$\Sing(h)\supset\Sing(V)$ the complement of the largest such Zariski
open set $X'$.
\endclaim

If $h$ is an admissible metric, we define 
$\cO_h(V^*)$ to be the sheaf of germs of holomorphic sections 
sections of $\smash{V^*_{|X\ssm\Sing(h)}}$ which are $h^*$-bounded 
near $\Sing(h)$; 
by the assumption on the analytic singularities, 
this is a coherent sheaf (as the direct image of some coherent sheaf on
$P(V)$), and actually, since $h^*=e^{-\varphi}h_0^*$, it is 
a subsheaf of the
sheaf $\cO(V^*):=\cO_{h_0}(V^*)$ associated with a smooth positive
definite metric $h_0$  on $T_X$. If $r$ is the generic rank of $V$ and
$m$ a positive integer, we define similarly $K_{V,h}^m$ to be sheaf of 
germs of holomorphic sections of 
$(\det V^*_{|X'})^{\otimes m}=(\Lambda^rV^*_{|X'})^{\otimes m}$ which 
are $\det h^*$-bounded, and $K_V^m:=K_{V,h_0}^m$.

If $V$ is defined by $\alpha:X\merto G_r(T_X)$, 
there always exists a modification $\mu:\smash{\wt X}\to X$ such that
the composition $\alpha\circ\mu:\smash{\wt X}\to G_r(\mu^* T_X)$
becomes holomorphic, and then $\mu^*V_{|\mu^{-1}(X\ssm \Sing(V))}$ 
extends as a locally trivial subbundle of $\mu^* T_X$ which we will
simply denote by $\mu^*V$. If $h$ is an admissible metric on $V$,
then $\mu^*V$ can be equipped with the
metric $\mu^*h=e^{\varphi\circ\mu}\mu^*h_0$ where $\mu^*h_0$
is smooth and positive definite. We may assume that $\varphi\circ\mu$
has divisorial singularities (otherwise just perform further blow-ups
of $\smash{\wt X}$ to achieve this). We then see that there is an integer 
$m_0$ such that for all multiples $m=pm_0$ the pull-back $\mu^*K_{V,h}^m$ 
is an invertible sheaf on $\smash{\wt X}$, and $\det h^*$ induces 
a smooth non singular metric on~it (when $h=h_0$, we can even take $m_0=1$).
By definition we always 
have $K_{V,h}^m=\mu_*(\mu^*K_{V,h}^m)$ for any $m\ge 0$.
In the sequel, however, we think of $K_{V,h}$ not really as a coherent 
sheaf, but rather as the ``virtual'' 
$\bQ$-line bundle $\mu_*\smash{(\mu^*K_{V,h}^{m_0})^{1/m_0}}$, and we say 
that $K_{V,h}$ is big
if $h^0(X, K_{V,h}^m)\ge cm^n$ for $m\ge m_1$, with $c>0\,$, i.e.\
if the invertible sheaf $\mu^*K_{V,h}^{m_0}$ is big in the usual sense.

At this point, it is important to observe that ``our'' canonical sheaf
$K_V$ differs from the sheaf $\cK_V:=i_*\cO(K_V)$ associated with
the injection $i:X\ssm \Sing(V)\hookrightarrow X$,
which is usually referred to as being the ``canonical sheaf'',
at least when $V$ is the space of tangents to a foliation. 
In fact, $\cK_V$ is always an invertible sheaf and there is
an obvious inclusion $K_V\subset\cK_V$. More precisely, the image of 
$\cO(\Lambda^rT^*_X)\to\cK_V$ is equal to $\cK_V\otimes_{\cO_X}\cJ$
for a certain coherent ideal $\cJ\subset\cO_X$, and the condition to 
have $h_0$-bounded sections on $X\ssm\Sing(V)$
precisely means that our sections are bounded by Const$\sum|g_j|$
in terms of the generators $(g_j)$ of $\cK_V\otimes_{\cO_X}\cJ$, i.e.\ 
$K_V=\cK_V\otimes_{\cO_X}\ol\cJ$ where $\ol\cJ$ is the integral
closure of~$\cJ$. More generally, 
$$
K_{V,h}^m=\cK_V^m\otimes_{\cO_X}\ol\cJ_{h,m_0}^{m/m_0}
$$
where $\ol\cJ_{h,m_0}^{m/m_0}\subset\cO_X$ is the $(m/m_0)$-integral closure 
of a certain ideal sheaf \hbox{$\cJ_{h,m_0}\subset\cO_X$},
which can itself be assumed to be integrally closed; in our \hbox{previous}
discussion, $\mu$ is chosen so that $\mu^*\cJ_{h,m_0}$ is invertible
on $\smash{\wt X}$.

The discrepancy already occurs
e.g.\ with the rank~$1$ linear space $V\subset T_{\bP^n_\bC}$ consisting at each
point $z\ne 0$ of the tangent to the line $(0z)$ (so that
necessarily $V_0=T_{\bP^n_\bC,0}$). As a sheaf (and not as a linear space), 
$i_*\cO(V)$ is the invertible sheaf generated by the vector field 
$\xi=\sum z_j{\partial/\partial z_j}$ on the affine open set
$\bC^n\subset\bP^n_\bC$, and therefore $\cK_V:=i_*\cO(V^*)$ is generated
over $\bC^n$ by the unique $1$-form $u$ such that $u(\xi)=1$. 
Since $\xi$ vanishes at $0$, the generator $u$ is {\it unbounded} 
with respect to a smooth metric $h_0$
on  $T_{\bP^n_\bC}$, and it is easily seen that $K_V$ is the non 
invertible sheaf $K_V=\cK_V\otimes \gm_{\bP^n_\bC,0}$. We can make it invertible
by considering the blow-up $\mu:\smash{\wt X}\to X$ of $X=\bP^n_\bC$
at $0$, so that $\mu^*K_V$ is isomorphic to $\mu^*\cK_V\otimes\cO_{\wt X}(-E)$ 
where $E$ is the exceptional divisor. The integral curves $C$ of $V$ are of
course lines through~$0$, and when a standard parametrization is used,
their derivatives do not vanish at~$0$, while the sections of 
$i_*\cO(V)$ do -- another sign
that $i_*\cO(V)$ and $i_*\cO(V^*)$ are the {\it wrong objects} to consider. 
Another standard example is obtained by taking a generic pencil of
elliptic curves $\lambda P(z)+\mu Q(z)=0$ of degree $3$ in $\bP_\bC^2$,
and the linear space $V$ consisting of the tangents to the fibers
of the rational map $\bP_\bC^2\merto\bP_\bC^1$
defined by $z\mapsto Q(z)/P(z)$. Then $V$ is given by
$$
0\lra i_*\cO(V)\longrightarrow\cO(T_{\bP_\bC^2})
\build\llra{8ex}^{PdQ-QdP}_{}\cO_{\bP^2_\bC}(6)\otimes\cJ_S\lra 0
$$
where $S=\Sing(V)$ consists of the 9 points $\{P(z)=0\}\cap\{Q(z)=0\}$, and
$\cJ_S$ is the corresponding ideal sheaf of~$S$. Since 
$\det\cO(T_{\bP^2})=\cO(3)$, we see that $\cK_V=\cO(3)$ is ample, which seems
to contradict (0.2) since all leaves are elliptic curves. There is however
no such contradiction, because $K_V=\cK_V\otimes\cJ_S$ is not big in our sense
(it has degree $0$ on all members of the elliptic pencil). A similar example
is obtained with a generic pencil of conics, in which case $\cK_V=\cO(1)$
and $\card S=4$.

For a given admissible hermitian structure $(V,h)$, we define similarly
the sheaf $E^\GG_{k, m}V^*_h$ to be the sheaf of 
polynomials defined over $X\ssm\Sing(h)$ which are \hbox{``$h$-bounded''}.
This means that when they are viewed as polynomials $P(z\,;\,\xi_1,\ldots,
\xi_k)$ in terms of $\xi_j=(\nabla_{h_0}^{1,0})^jf(0)$ where
$\nabla_{h_0}^{1,0}$ is the $(1,0)$-component of the induced Chern 
connection on $(V,h_0)$, there is a uniform bound 
$$
\big|P(z\,;\,\xi_1,\ldots,\xi_k)\big|\le 
C\Big(\sum\Vert \xi_j\Vert_h^{1/j}\Big)^m\leqno(0.12)
$$
near points of $X\ssm X'$ (see section~2 for more details on
this). Again, by a direct image argument, one sees that
$\smash{E^\GG_{k, m}V^*_h}$ is always a coherent sheaf. The sheaf
$\smash{E^\GG_{k, m}V^*}$ is defined to be $\smash{E^\GG_{k, m}V^*_h}$ 
when $h=h_0$ (it is actually independent of the choice of $h_0$, as
follows from arguments similar to those given in section~2). Notice
that this is exactly what is needed to extend the proof of the vanishing
theorem~0.4 to the case of a singular linear space $V\,$; the value 
distribution theory argument can only work when the functions 
$P(f\,;\;f',\ldots,f^{(k)})(t)$ do not exhibit poles, and this is guaranteed
here by the boundedness assumption.

Our strategy can be described as follows. We consider the
Green-Griffiths bundle of $k$-jets $X_k^\GG=J^kV\ssm\{0\}/\bC^*$,
which by (0.3) consists of a fibration in {\it weighted projective
spaces}, and its associated tautological sheaf 
$$L=\cO_{X_k^\GG}(1),$$
viewed rather as a virtual $\bQ$-line bundle $\smash{\cO_{X_k^\GG}}
(m_0)^{1/m_0}$ with $m_0=\lcm(1,2,\,...\,,k)$.
Then, if $\pi_k:\smash{X_k^\GG}\to X$ is the natural projection, we have
$$
E^\GG_{k,m}=(\pi_k)_*\cO_{X_k^\GG}(m)\quad\hbox{and}\quad
R^q(\pi_k)_*\cO_{X_k^\GG}(m)=0~\hbox{for $q\ge 1$}.
$$
Hence, by the Leray spectral sequence we get for every invertible sheaf
$F$ on $X$ the isomorphism
$$
H^q(X,E^\GG_{k,m}V^*\otimes F)\simeq 
H^q(X_k^\GG,\cO_{X_k^\GG}(m)\otimes\pi_k^*F).
\leqno(0.13)
$$
The latter group can be evaluated thanks to holomorphic Morse inequalities.
In fact we can associate with any admissible metric $h$ on $V$ a metric
(or rather a natural family) of metrics on $L=\smash{\cO_{X_k^\GG}}(1)$. 
The space $\smash{X_k^\GG}$ always possesses quotient singularities if $k\ge 2$ 
(and even some more if $V$ is singular), but we do not really care since 
Morse inequalities still work in this setting. As we will see, it is then
possible to get nice asymptotic formulas as $k\to +\infty$. They appear to
be of a {\it probabilistic nature} if we take the components of the $k$-jet
(i.e.\ the successive derivatives $\xi_j=f^{(j)}(0)$, $1\le j\le k$) as 
random variables. This probabilistic behaviour was somehow  already
visible in the Riemann-Roch calculation of~[GG79].
In this way, assuming $K_V$ big, we produce a lot of sections
$\sigma_j=\smash{H^0(X_k^\GG,\cO_{X_k^\GG}}(m)\otimes\pi_k^*F)$, corresponding
to certain divisors $Z_j\subset \smash{X_k^\GG}$. The hard problem which
is left in order to complete a proof of the generalized
Green-Griffiths-Lang conjecture is to compute the base locus
$Z=\bigcap Z_j$ and to show that $Y=\pi_k(Z)\subset X$ must be a
proper algebraic variety. Although we cannot address this problem at
present, we will indicate a few technical results and a couple of
potential strategies in this direction.

I would like to thank Simone Diverio and Mihai Paun for several stimulating
discussions, and Erwan Rousseau for convincing me to explain better
the peculiarities of the definition of the canonical sheaf employed here.

\section{1. Hermitian geometry of weighted projective spaces}

The goal of this section is to introduce natural K\"ahler metrics on
weighted projective spaces, and to evaluate the corresponding volume forms.
Here we put $d^c={i\over 4\pi}(\dbar-\partial)$ so that
$dd^c={i\over 2\pi}\ddbar$. The normalization of the $d^c$ operator is
chosen such that we have precisely $(dd^c\log|z|^2)^n=\delta_0$ for the
Monge-Amp\`ere operator in $\bC^n$; also, for every holomorphic or meromorphic
section $\sigma$ of a hermitian line bundle $(L,h)$  the Lelong-Poincar\'e
can be formulated
$$dd^c\log|\sigma|^2_h=[Z_\sigma]-\Theta_{L,h},\leqno(1.1)$$ 
where $\Theta_{L,h}={i\over 2\pi}D_{L,h}^2$ is the $(1,1)$-curvature 
form of $L$ and $Z_\sigma$ 
the zero divisor of~$\sigma$.  The closed $(1,1)$-form $\Theta_{L,h}$ is
a representative of the first Chern class $c_1(L)$. Given a $k$-tuple
of ``weights'' $a=(a_1,\ldots,a_k)$, i.e.\ of relatively prime
integers $a_s>0$ we introduce the weighted projective space
$P(a_1,\ldots,a_k)$ to be the quotient of $\bC^k\ssm\{0\}$ by the
corresponding weighted $\bC^*$ action:
$$
P(a_1,\ldots,a_k)=\bC^k\ssm\{0\}/\bC^*,\qquad
\lambda\cdot z =(\lambda^{a_1}z_1,\ldots,\lambda^{a_k}z_k).
\leqno(1.2)
$$
As is well known, this defines a toric $k-1$-dimensional algebraic variety 
with quotient singularities. On this variety, we introduce the possibly
singular (but almost everywhere smooth and non degenerate) K\"ahler form 
$\omega_{a,p}$ defined by
$$
\pi_a^*\omega_{a,p}=dd^c\varphi_{a,p},\qquad
\varphi_{a,p}(z)={1\over p}\log\sum_{1\le s\le k}|z_s|^{2p/a_s},
\leqno(1.3)
$$
where $\pi_a:\bC^k\ssm\{0\}\to P(a_1,\ldots,a_k)$ is the canonical projection
and $p>0$ is a positive constant. It is clear that $\varphi_{p,a}$ is real
analytic on $\bC^k\ssm\{0\}$ if $p$ is an integer and a common multiple of 
all weights $a_s$. It is at least $C^2$ is $p$ is real and $p\ge\max(a_s)$,
which will be more than sufficient for our purposes (but everything
would still work for any $p>0$). The resulting metric is in any case 
smooth and positive definite outside of the coordinate hyperplanes
$z_s=0$, and these hyperplanes will not matter here since they are of 
capacity zero with respect to all currents $(dd^c\varphi_{a,p})^\ell$.
In order to evaluate the volume $\int_{P(a_1,\ldots,a_k)}\omega_{a,p}^{k-1}$,
one can observe that
$$
\leqalignno{
\int_{P(a_1,\ldots,a_k)}\omega_{a,p}^{k-1}
&=\int_{z\in\bC^k,\,\varphi_{a,p}(z)=0}\pi_a^*\omega_{a,p}^{k-1}\wedge
d^c\varphi_{a,p}\cr
&=
\int_{z\in\bC^k,\,\varphi_{a,p}(z)=0}(dd^c\varphi_{a,p})^{k-1}\wedge
d^c\varphi_{a,p}\cr
&={1\over p^k}\int_{z\in\bC^k,\,\varphi_{a,p}(z)<0}(dd^ce^{p\varphi_{a,p}})^k.
&(1.4)\cr}
$$
The first equality comes from the fact that $\{\varphi_{a,p}(z)=0\}$ is a
circle bundle over $P(a_1,\ldots,a_k)$, by using the identities
$\varphi_{a,p}(\lambda\cdot z)=\varphi_{a,p}(z)+\log|\lambda|^2$
and $\int_{|\lambda|=1}d^c\log|\lambda|^2=1$. The third equality can be seen
by Stokes formula applied to the $(2k-1)$-form
$$
(dd^ce^{p\varphi_{a,p}})^{k-1}\wedge
d^ce^{p\varphi_{a,p}}=
e^{p\varphi_{a,p}}(dd^c\varphi_{a,p})^{k-1}\wedge d^c\varphi_{a,p}
$$
on the pseudoconvex open set $\{z\in\bC^k\,;\;\varphi_{a,p}(z)<0\}$. 
Now, we find
$$
\leqalignno{\kern20pt
(dd^ce^{p\varphi_{a,p}})^k=\Big(dd^c\sum_{1\le s\le k}|z_s|^{2p/a_s}\Big)^k
&=\prod_{1\le s\le k}
\Big({p\over a_s}|z_s|^{{p\over a_s}-1}\Big)(dd^c|z|^2)^k,&(1.5)\cr
\int_{z\in\bC^k,\,\varphi_{a,p}(z)<0}
(dd^ce^{p\varphi_{a,p}})^k&=\prod_{1\le s\le k}{p\over a_s}=
{p^k\over a_1\ldots a_k}.&(1.6)\cr}
$$
In fact, (1.5) and (1.6) are clear when $p=a_1=\ldots=a_k=1$ (this is
just the standard calculation of the volume of the unit ball in $\bC^k$);
the general case follows by substituting formally $z_s\mapsto 
\smash{z_s^{p/a_s}}$, and using rotational invariance along with the
observation that the arguments of the complex
numbers $\smash{z_s^{p/a_s}}$ now run in the interval $[0,2\pi p/a_s[$ 
instead of $[0,2\pi[$ (say). As a consequence of
(1.4) and (1.6), we 
obtain the well known value
$$
\int_{P(a_1,\ldots,a_k)}\omega_{a,p}^{k-1}={1\over a_1\ldots a_k},\leqno(1.7)
$$
for the volume. Notice that this is independent of $p$ (as it is obvious 
by Stokes theorem, since the cohomology class of $\omega_{a,p}$ does 
not depend on $p$).
When $p$ tends to $+\infty$, we have $\varphi_{a,p}(z)\mapsto
\varphi_{a,\infty}(z)=\log\max_{1\le s\le k}|z_s|^{2/a_s}$ and the volume form 
$\omega_{a,p}^{k-1}$ converges to a rotationally invariant measure 
supported by the image of the polycircle $\prod\{|z_s|=1\}$ in
$P(a_1,\ldots,a_k)$. This is so because not all $|z_s|^{2/a_s}$ are equal
outside of the image of the polycircle, thus $\varphi_{a,\infty}(z)$
locally depends only on~$k-1$ complex variables, and so 
$\omega_{a,\infty}^{k-1}=0$ there by log homogeneity.

Our later calculations will require a slightly more general setting.
Instead of looking at $\bC^k$, we consider the weighted $\bC^*$ action
defined by
$$
\bC^{|r|}=\bC^{r_1}\times\ldots\times\bC^{r_k},\qquad
\lambda\cdot z =(\lambda^{a_1}z_1,\ldots,\lambda^{a_k}z_k).
\leqno(1.8)
$$
Here $z_s\in\bC^{r_s}$ for some $k$-tuple $r=(r_1,\ldots,r_k)$ and
$|r|=r_1+\ldots+r_k$. This gives rise to a weighted projective space 
$$
\leqalignno{
&P(a_1^{[r_1]},\ldots,a_k^{[r_k]})=P(a_1,\ldots,a_1,\ldots,a_k,\ldots,a_k),\cr
&\pi_{a,r}:\bC^{r_1}\times\ldots\times\bC^{r_k}
\ssm\{0\}\longrightarrow P(a_1^{[r_1]},\ldots,a_k^{[r_k]})&(1.9)\cr}
$$
obtained by repeating $r_s$ times each weight $a_s$. On this space, we
introduce the degenerate K\"ahler metric $\omega_{a,r,p}$ such that
$$
\pi_{a,r}^*\omega_{a,r,p}=dd^c\varphi_{a,r,p},\qquad
\varphi_{a,r,p}(z)={1\over p}\log\sum_{1\le s\le k}|z_s|^{2p/a_s}
\leqno(1.10)
$$
where $|z_s|$ stands now for the standard hermitian norm 
$(\sum_{1\le j\le r_s}|z_{s,j}|^2)^{1/2}$ on~$\bC^{r_s}$. This metric
is cohomologous to the corresponding ``polydisc-like'' 
metric $\omega_{a,p}$ already defined, and therefore Stokes theorem implies
$$
\int_{P(a_1^{[r_1]},\ldots,a_k^{[r_k]})}
\omega_{a,r,p}^{|r|-1}={1\over a_1^{r_1}\ldots a_k^{r_k}}.
\leqno(1.11)
$$
Since $(dd^c\log|z_s|^2)^{r_s}=0$ on $\bC^{r_s}\ssm\{0\}$ by homogeneity,
we conclude as before that the weak limit
$\lim_{p\to+\infty}\omega_{a,r,p}^{|r|-1}=\omega_{a,r,\infty}^{|r|-1}$ associated
with
$$
\varphi_{a,r,\infty}(z)=\log\max_{1\le s\le k}|z_s|^{2/a_s}\leqno(1.12)
$$
is a measure supported by the image of the product of unit spheres
$\prod S^{2r_s-1}$ in $\smash{P(a_1^{[r_1]},\ldots,a_k^{[r_k]})}$, which is
invariant under the action of 
\hbox{$U(r_1)\times\ldots\times U(r_k)$} on 
\hbox{$\bC^{r_1}\times\ldots\times\bC^{r_k}$},
and thus coincides with the hermitian area measure up to a constant 
determined by condition (1.11). In fact, outside of the product of
spheres, $\varphi_{a,r,\infty}$ locally depends only on at most
$k-1$ factors and thus, for dimension reasons, the top power 
$(dd^c\varphi_{a,r,\infty})^{|r|-1}$ must be zero there. In the next section,
the following change of variable formula will be needed. For simplicity
of exposition we restrict ourselves to continuous functions, but 
a standard density argument would easily extend the formula to all
functions that are Lebesgue integrable with respect to the volume 
form~$\smash{\omega_{a,r,p}^{|r|-1}}$.

\claim (1.13) Proposition|Let $f(z)$ be a bounded function on
$P(a_1^{[r_1]},\ldots,a_k^{[r_k]})$ which is continuous outside of
the hyperplane sections $z_s=0$. We also view $f$ as a $\bC^*$-invariant
continuous function on $\prod(\bC^{r_s}\ssm\{0\})$. Then
$$
\eqalign{
&\int_{P(a_1^{[r_1]},\ldots,a_k^{[r_k]})}f(z)\,\omega_{a,r,p}^{|r|-1}\cr
&={(|r|-1)!\over \prod_s a_s^{r_s}}\kern-2pt
\int_{(x,u)\in\Delta_{k-1}\times\prod S^{2r_s-1}}\kern-3pt
f(x_1^{a_1/2p}u_1,\ldots,x_k^{a_k/2p}u_k)\kern-3.5pt
\prod_{1\le s\le k}\kern-2.5pt
{x_s^{r_s-1}\over(r_s-1)!}\,dx\,d\mu(u)\kern-3pt
\cr}
$$
where $\Delta_{k-1}$ is the $(k-1)$-simplex $\{x_s\ge 0$, $\sum x_s=1\}$,
$dx=dx_1\wedge\ldots\wedge dx_{k-1}$ its standard measure, and where
$d\mu(u)=d\mu_1(u_1)\ldots d\mu_k(u_k)$ is the rotation invariant 
probability measure on the product $\prod_sS^{2r_s-1}$ of unit 
spheres in $\bC^{r_1}\times\ldots\times\bC^{r_k}$. As a consequence
$$
\lim_{p\to+\infty}
\int_{P(a_1^{[r_1]},\ldots,a_k^{[r_k]})}f(z)\,\omega_{a,r,p}^{|r|-1}=
{1\over \prod_s a_s^{r_s}}\int_{\prod S^{2r_s-1}}f(u)\,d\mu(u).
$$
\endclaim

\proof. The area formula of the disc $\smash{\int_{|\lambda|<1}}
dd^c|\lambda|^2=1$ and a consideration of the unit disc bundle over
$\smash{P(a_1^{[r_1]},\ldots,a_k^{[r_k]})}$ imply that
$$
I_p:=\int_{P(a_1^{[r_1]},\ldots,a_k^{[r_k]})}\kern-2pt
f(z)\,\omega_{a,r,p}^{|r|-1}=\int_{z\in\bC^{|r|},\varphi_{a,r,p}(z)<0}
\kern-2pt f(z)\,
(dd^c\varphi_{a,r,p})^{|r|-1}\wedge dd^ce^{\varphi_{a,r,p}}.\kern-2pt
$$
Now, a straightforward calculation on $\bC^{|r|}$ gives
$$
\eqalign{
(dd^ce^{p\varphi_{a,r,p}})^{|r|}
&=\Big(dd^c\sum_{1\le s\le k}|z_s|^{2p/a_s}\Big)^{|r|}\cr
&=\prod_{1\le s\le k}\Big({p\over a_s}\Big)^{r_s+1}|z_s|^{2r_s(p/a_s-1)}
(dd^c|z|^2)^{|r|}.\cr}
$$
On the other hand, we have $(dd^c|z|^2)^{|r|}={|r|!\over r_1!\ldots r_k!}
\prod_{1\le s\le k}(dd^c|z_s|^2)^{r_s}$ and
$$
\leqalignno{
(dd^ce^{p\varphi_{a,r,p}})^{|r|}
&=\big(p\,e^{p\varphi_{a,r,p}}(dd^c\varphi_{a,r,p}+
p\,d\varphi_{a,r,p}\wedge d^c\varphi_{a,r,p})\big)^{|r|}\cr
&=|r|p^{|r|+1}e^{|r|p\varphi_{a,r,p}}(dd^c\varphi_{a,r,p})^{|r|-1}\wedge
d\varphi_{a,r,p}\wedge d^c\varphi_{a,r,p}\cr
&=|r|p^{|r|+1}e^{(|r|p-1)\varphi_{a,r,p}}(dd^c\varphi_{a,r,p})^{|r|-1}\wedge
dd^ce^{\varphi_{a,r,p}},\cr}
$$
thanks to the homogeneity relation $(dd^c\varphi_{a,r,p})^{|r|}=0$. Putting 
everything together, we find
$$
I_p=\int_{z\in\bC^{|r|},\,\varphi_{a,r,p}(z)<0}~~
{(|r|-1)!\,p^{k-1}f(z)\over (\sum_s|z_s|^{2p/a_s})^{|r|-1/p}}
\,\prod_s{(dd^c|z_s|^2)^{r_s}\over r_s!\,a_s^{r_s+1}|z_s|^{2r_s(1-p/a_s)}}.
$$
A standard calculation in polar coordinates with $z_s=\rho_su_s$,
$u_s\in S^{2r_s-1}$, yields
$$
{(dd^c|z_s|^2)^{r_s}\over |z_s|^{2r_s}}=2r_s{d\rho_s\over \rho_s}\,d\mu_s(u_s)
$$
where $\mu_s$ is the $U(r_s)$-invariant probability measure on 
$S^{2r_s-1}$. Therefore
$$
\eqalign{I_p
&=\int_{\varphi_{a,r,p}(z)<0}
{(|r|-1)!\,p^{k-1}f(\rho_1u_1,\ldots,\rho_ku_k)\over
(\sum_{1\le s\le k}\rho_s^{2p/a_s})^{|r|-1/p}}
\prod_s{2\rho_s^{2pr_s/a_s}{d\rho_s\over \rho_s}d\mu_s(u_s)
\over (r_s-1)!\,a_s^{r_s+1}}\cr
&=\int_{u_s\in S^{2r_s-1},\,\sum t_s<1}\kern-3pt
{(|r|-1)!\,p^{-1}f(t_1^{a_1/2p}u_1,\ldots,t_k^{a_k/2p}u_k)\over
(\sum_{1\le s\le k}t_s)^{|r|-1/p}}
\prod_s{t_s^{r_s-1}dt_s\,d\mu_s(u_s)
\over (r_s-1)!\,a_s^{r_s}}\cr}
$$
by putting $t_s=|z_s|^{2p/a_s}=\rho_s^{2p/a_s}$, i.e.\ 
$\rho_s=t_s^{a_s/2p}$, $t_s\in{}]0,1]$. We use still another change of variable 
$t_s=tx_s$ with $t=\sum_{1\le s\le k}t_s$ and $x_s\in{}]0,1]$,
$\sum_{1\le s\le k}x_s=1$. Then
$$dt_1\wedge\ldots\wedge dt_k=t^{k-1}\,dx\,dt\qquad
\hbox{where $dx=dx_1\wedge\ldots\wedge dx_{k-1}$}.
$$
The $\bC^*$ invariance of $f$ shows that
$$
\eqalign{
I_p&=\int_{u_s\in S^{2r_s-1}\atop{\scriptstyle \Sigma} x_s=1,\, t\in{}]0,1]}
(|r|-1)!f(x_1^{a_s/2p}u_1,\ldots,x_k^{a_k/2p}u_k)
\prod_{1\le s\le k}{x_s^{r_s-1}d\mu_s(u_s)\over (r_s-1)!\,a_s^{r_s}}
{dx\,dt\over p\,t^{1-1/p}}\cr
&=\int_{u_s\in S^{2r_s-1}\atop{\scriptstyle\Sigma}x_s=1\kern18pt}
(|r|-1)!f(x_1^{a_s/2p}u_1,\ldots,x_k^{a_k/2p}u_k)
\prod_{1\le s\le k}{x_s^{r_s-1}d\mu_s(u_s)\over (r_s-1)!\,a_s^{r_s}}\,dx.\cr}
$$
This is equivalent to the formula given in Proposition 1.13. We have
$x_s^{2a_s/p}\to 1$ as $p\to+\infty$, and by Lebesgue's bounded convergence 
theorem and Fubini's formula, we get
$$
\lim_{p\to+\infty}I_p=
{(|r|-1)!\over\prod_s a_s^{r_s}}
\int_{(x,u)\in\Delta_{k-1}\times\prod S^{2r_s-1}}
f(u)\prod_{1\le s\le k}{x_s^{r_s-1}\over(r_s-1)!}\,dx\,d\mu(u).
$$
It can be checked by elementary integrations by parts and induction on 
$k,\,r_1,\ldots,r_k$ that
$$
\int_{x\in\Delta_{k-1}}
\prod_{1\le s\le k}x_s^{r_s-1}dx_1\ldots dx_{k-1}
={1\over (|r|-1)!}\prod_{1\le s\le k}(r_s-1)!~.
\leqno(1.14)
$$
This implies that $(|r|-1)!\prod_{1\le s\le k}
{x_s^{r_s-1}\over(r_s-1)!}\,dx$ is a probability measure on $\Delta_{k-1}$ 
and that
$$
\lim_{p\to+\infty}I_p=
{1\over\prod_s a_s^{r_s}}\int_{u\in\prod S^{2r_s-1}}f(u)\,d\mu(u).
$$
Even without an explicit check, the evaluation (1.14) also follows
from the fact that we must have equality for $f(z)\equiv 1$ in the
latter equality, if we take into account the volume formula
(1.11).\qed

\section{2. Probabilistic estimate of the curvature of 
{\twelvebfit k}-jet bundles}
\sectionrunning{2. Probabilistic estimate of the curvature of $k$-jet bundles}

Let $(X,V)$ be a compact complex directed non singular variety. To
avoid any technical difficulty at this point, we first assume that $V$
is a holomorphic vector subbundle of $T_X$, equipped with a smooth
hermitian metric $h$.

According to the notation already specified in the introduction, we 
denote by $J^kV$ the bundle of $k$-jets of holomorphic curves
$f:(\bC,0)\to X$ tangent to $V$ at each point. Let us set $n=\dim_\bC
X$ and $r=\rank_\bC V$. Then $J^kV\to X$ is an algebraic fiber bundle
with typical fiber $\bC^{rk}$ (see below). It has a canonical
$\bC^*$-action defined by $\lambda\cdot f:(\bC,0)\to X$,
$(\lambda\cdot f)(t)=f(\lambda t)$.  Fix a point $x_0$ in $X$ and a
local holomorphic coordinate system $(z_1,\ldots,z_n)$ centered at
$x_0$ such that $V_{x_0}$ is the vector subspace
$\langle \partial/\partial z_1,\ldots,
\partial/\partial z_r\rangle$ at $x_0$. Then, in a neighborhood $U$ of $x_0$,
$V$ admits a holomorphic frame of the form
$$
{\partial\over \partial z_\beta}+\sum_{r+1\le\alpha\le n}
a_{\alpha\beta}(z){\partial\over \partial z_\alpha},\qquad
1\le\beta\le r,~~a_{\alpha\beta}(0)=0.
\leqno(2.1)
$$
Let $f(t)=(f_1(t),\ldots,f_n(t))$ be a $k$-jet of curve tangent to $V$ starting 
from a point $f(0)=x\in U$. Such a curve is entirely determined by its 
initial point and by the projection $\wt f(t):=(f_1(t),\ldots,f_r(t))$ 
to the first
$r$-components, since the condition $f'(t)\in V_{f(t)}$ implies that the
other components must satisfy the ordinary differential equation
$$
f'_\alpha(t)=\sum_{1\le\beta\le r}a_{\alpha\beta}(f(t))f'_\beta(t).
$$
This implies that the $k$-jet of $f$ is entirely determined by 
the initial point $x$ and the Taylor expansion
$$
\wt f(t)-\wt x=\xi_1t+\xi_2t^2+\ldots+\xi_kt^k+O(t^{k+1})
\leqno(2.2)
$$
where $\xi_s=(\xi_{s\alpha})_{1\le\alpha\le r}\in\bC^r$. The $\bC^*$ action
$(\lambda,f)\mapsto \lambda\cdot f$ is then expressed in coordinates
by the weighted action
$$
\lambda\cdot(\xi_1,\xi_2,\ldots,\xi_k)=
(\lambda\xi_1,\lambda^2\xi_2,\ldots,\lambda^k\xi_k)
\leqno(2.3)
$$
associated with the weight $a=(1^{[r]},2^{[r]},\ldots,k^{[r]})$. The quotient
projectived $k$-jet bundle 
$$
X^\GG_k:=(J^kV\ssm\{0\})/\bC^*
\leqno(2.4)
$$
considered by Green and Griffiths [GG79] is therefore in a natural
way a $P(1^{[r]},2^{[r]},\ldots,k^{[r]})$ weighted projective bundle over $X$.
As such, it possesses a canonical sheaf $\smash{\cO_{X^\GG_k}}(1)$ such that
$\smash{\cO_{X^\GG_k}}(m)$ is invertible when $m$ is a multiple of
$\lcm(1,2,\ldots,k)$. Under the natural projection 
$\pi_k:\smash{X^\GG_k}\to X$, the direct image 
$(\pi_k)_*\smash{\cO_{X^\GG_k}(m)}$ coincides with the 
sheaf of sections of the bundle $E_{k,m}^\GG V^*$ of jet differentials 
of order $k$ and degree $m$, namely polynomials
$$
P(z\,;\,\xi_1,\ldots,\xi_k)=\sum_{\alpha_\ell\in\bN^r,\,1\le\ell\le k} 
a_{\alpha_1\ldots\alpha_k}(z)\,\xi_1^{\alpha_1}\ldots\xi_k^{\alpha_k}
\leqno(2.5)
$$
of weighted degree $|\alpha_1|+2|\alpha_2|+\ldots+k|\alpha_k|=m$ on 
$J^kV$ with holomorphic coefficients. The jet differentials operate
on germs of curves as differential operators
$$
P(f)(t)=\sum a_{\alpha_1\ldots\alpha_k}(f(t))\;
f'(t)^{\alpha_1}\ldots f^{(k)}(t)^{\alpha_k}
\leqno(2.6)
$$
In the sequel, we do not make any further use of coordinate frames as
(2.1), because they need not be related in any way to the hermitian 
metric $h$ of $V$.  Instead, we choose a local holomorphic coordinate frame
$(e_\alpha(z))_{1\le\alpha\le r}$ of $V$ on a neighborhood $U$ of~$x_0$, 
such that
$$
\langle e_\alpha(z),e_\beta(z)\rangle =\delta_{\alpha\beta}+
\sum_{1\le i,j\le n,\,1\le\alpha,\beta\le r}c_{ij\alpha\beta}z_i\overline z_j+
O(|z|^3)\leqno(2.7)
$$
for suitable complex coefficients $(c_{ij\alpha\beta})$. It is a standard fact
that such a normalized coordinate system always exists, and that the 
Chern curvature tensor ${i\over 2\pi}D^2_{V,h}$ of $(V,h)$ at $x_0$ 
is then given by
$$
\Theta_{V,h}(x_0)=-{i\over 2\pi}
\sum_{i,j,\alpha,\beta}
c_{ij\alpha\beta}\,dz_i\wedge d\overline z_j\otimes e_\alpha^*\otimes e_\beta.
\leqno(2.8)
$$
Also, instead of defining the vectors $\xi_s\in\bC^r$ as in (2.2),
we consider a local holomorphic connection $\nabla$ on $V_{|U}$ (e.g.\ the 
one which turns $(e_\alpha)$ into a parallel frame), and take 
$\xi_k=\nabla^kf(0)\in V_x$ defined inductively 
by $\nabla^1 f=f'$ and $\nabla^sf=\nabla_{f'}(\nabla^{s-1}f)$. This is
just another way of parametrizing the fibers of $J^kV$ over $U$ by the
vector bundle $\smash{V^k_{|U}}$. Notice that this is highly dependent 
on $\nabla$
(the bundle $J^kV$ actually does not carry a vector bundle or even affine
bundle structure); however, the expression of the weighted action (2.3) 
is unchanged in this new setting. Now, we fix a finite open covering 
$(U_\alpha)_{\alpha\in I}$ of~$X$ by open coordinated charts such that
$V_{|U_\alpha}$ is trivial, along with holomorphic connections 
$\nabla_\alpha$ on $V_{|U_\alpha}$. Let $\theta_\alpha$ be a partition of
unity of $X$ subordinate to the covering $(U_\alpha)$. Let us fix 
$p>0$ and small parameters $1=\varepsilon_1\gg\varepsilon_2\gg\ldots\gg
\varepsilon_k>0$. Then we define a global 
weighted exhaustion on $J^kV$ by putting for any $k$-jet $f\in J^k_xV$
$$
\Psi_{h,p,\varepsilon}(f):=\Big(\sum_{\alpha\in I}
\theta_\alpha(x)\sum_{1\le s\le k}\varepsilon_s^{2p}\Vert\nabla^s_\alpha f(0)
\Vert_{h(x)}^{2p/s}\Big)^{1/p}
\leqno(2.9)
$$
where $\Vert~~\Vert_{h(x)}$ is the hermitian metric $h$ of $V$ evaluated
on the fiber $V_x$, $x=f(0)$. The function $\Psi_{h,p,\varepsilon}$ satisfies
the fundamental homogeneity property 
$$
\Psi_{h,p,\varepsilon}(\lambda\cdot f)=\Psi_{h,p,\varepsilon}(f)\,|\lambda|^2
\leqno(2.10)
$$
with respect to the $\bC^*$ action on $J^kV$, in other words, it induces
a hermitian metric on the dual $L^*$ of the tautological $\bQ$-line bundle
$L_k=\cO_{X_k^\GG}(1)$ over $X_k^\GG$. The curvature of $L_k$ is given by
$$
\pi_k^*\Theta_{L_k,\Psi^*_{h,p,\varepsilon}}=dd^c\log\Psi_{h,p,\varepsilon}
\leqno(2.11)
$$
where $\pi_k:J^kV\ssm\{0\}\to X_k^\GG$ is the canonical projection. 
Our next goal is to compute precisely the curvature and to apply
holomorphic Morse inequalities to $L\to X_k^\GG$ with the above metric.
It might look a priori like an untractable problem, since the definition of
$\Psi_{h,p,\varepsilon}$ is a rather unnatural one. However, the ``miracle''
is that the asymptotic behavior of $\Psi_{h,p,\varepsilon}$ as
$\varepsilon_s/\varepsilon_{s-1}\to 0$ is in some sense uniquely defined 
and very natural.
It will lead to a computable asymptotic formula, which is moreover
simple enough to produce useful results.

\vbox{%
\claim (2.12) Lemma|On each coordinate chart $U$ equipped with
a holomorphic connection $\nabla$ of $V_{|U}$, let us define 
the components of a $k$-jet $f\in J^kV$ by $\xi_s=\nabla^sf(0)$,
and consider the rescaling transformation 
$$\rho_\varepsilon(\xi_1,\xi_2,\ldots,\xi_k)=
(\varepsilon_1^1\xi_1,\varepsilon_2^2\xi_2,\ldots,
\varepsilon_k^k\xi_k)\quad
\hbox{on $J^k_xV$, $x\in U_\alpha$}
$$
$($it commutes with the $\bC^*$-action but is otherwise unrelated and 
not canonically defined over $X$ as it depends on the choice of 
$\nabla)$. Then, if $p$ is a multiple of $\lcm(1,2,\ldots,k)$ and
$\varepsilon_s/\varepsilon_{s-1}\to 0$ for all $s=2,\ldots,k$, the
rescaled function $\Psi_{h,p,\varepsilon}\circ\rho_\varepsilon^{-1}
(\xi_1,\ldots,\xi_k)$ converges towards
$$
\smash{\bigg(\sum_{1\le s\le k}\Vert \xi_s\Vert^{2p/s}_h\bigg)^{1/p}}
$$
on every compact subset of $J^kV_{|U}\ssm\{0\}$,
uniformly in $C^\infty$ topology. \dummy
\endclaim\vskip-15pt}

\proof. Let $U\subset X$ be an open set on which $V_{|U}$ is trivial
and equipped with some holomorphic connection $\nabla$. Let us pick
another holomorphic connection $\wt\nabla=
\nabla+\Gamma$ where $\Gamma\in H^0(U,\Omega^1_X\otimes
\Hom(V,V)$. Then $\wt\nabla^2f=\nabla^2f+\Gamma(f)(f')\cdot f'$, and
inductively we get
$$
\wt\nabla^sf=\nabla^sf+P_s(f\,;\,\nabla^1f,\ldots,\nabla^{s-1}f)
$$
where $P(x\,;\,\xi_1,\ldots,\xi_{s-1})$ is a polynomial with holomorphic 
coefficients in $x\in U$ which is of weighted homogeneous degree
$s$ in $(\xi_1,\ldots,\xi_{s-1})$. In other words, the corresponding change
in the parametrization of $J^kV_{|U}$ is given by a $\bC^*$-homogeneous
transformation
$$
\wt\xi_s=\xi_s+P_s(x\,;\,\xi_1,\ldots,\xi_{s-1}).
$$
Let us introduce the corresponding rescaled components
$$
(\xi_{1,\varepsilon},\ldots,\xi_{k,\varepsilon})=
(\varepsilon_1^1\xi_1,\ldots,\varepsilon_k^k\xi_k),\qquad
(\wt\xi_{1,\varepsilon},\ldots,\wt\xi_{k,\varepsilon})=
(\varepsilon_1^1\wt\xi_1,\ldots,\varepsilon_k^k\wt\xi_k).
$$
Then
$$
\eqalign{
\wt\xi_{s,\varepsilon}
&=\xi_{s,\varepsilon}+
\varepsilon_s^s\,P_s(x\,;\,\varepsilon_1^{-1}\xi_{1,\varepsilon},\ldots,
\varepsilon_{s-1}^{-(s-1)}\xi_{s-1,\varepsilon})\cr
&=\xi_{s,\varepsilon}+O(\varepsilon_s/\varepsilon_{s-1})^s\,
O(\Vert\xi_{1,\varepsilon}\Vert+\ldots+\Vert\xi_{s-1,\varepsilon}
\Vert^{1/(s-1)})^s\cr}
$$
and the error terms are thus polynomials of fixed degree with arbitrarily
small coefficients as $\varepsilon_s/\varepsilon_{s-1}\to 0$. Now, the 
definition of $\Psi_{h,p,\varepsilon}$ consists of glueing the sums
$$
\sum_{1\le s\le k}\varepsilon_s^{2p}\Vert\xi_k\Vert_h^{2p/s}=
\sum_{1\le s\le k}\Vert\xi_{k,\varepsilon}\Vert_h^{2p/s}
$$
corresponding to $\xi_k=\nabla_\alpha^sf(0)$ by means of the partition
of unity $\sum\theta_\alpha(x)=1$. We see that by using the rescaled
variables $\xi_{s,\varepsilon}$ the changes occurring when replacing a
connection $\nabla_\alpha$ by an alternative one $\nabla_\beta$ are
arbitrary small in $C^\infty$ topology, with error terms uniformly
controlled in terms of the ratios $\varepsilon_s/\varepsilon_{s-1}$ on
all compact subsets of $V^k\ssm\{0\}$. This shows that in $C^\infty$
topology,
$\Psi_{h,p,\varepsilon}\circ\rho_\varepsilon^{-1}(\xi_1,\ldots,\xi_k)$
converges uniformly towards $\smash{(\sum_{1\le s\le
    k}\Vert\xi_k\Vert_h^{2p/s})^{1/p}}$, whatever is the trivializing
open set $U$ and the holomorphic connection $\nabla$ used to evaluate
the components and perform the rescaling.\qed

Now, we fix a point $x_0\in X$ and a local holomorphic frame 
$(e_\alpha(z))_{1\le\alpha\le r}$ satisfying (2.7) on a neighborhood $U$ 
of~$x_0$. We introduce the rescaled components 
$\xi_s=\varepsilon_s^s\nabla^sf(0)$ on $J^kV_{|U}$ and compute
the curvature of
$$
\Psi_{h,p,\varepsilon}\circ\rho_\varepsilon^{-1}(z\,;\,\xi_1,\ldots,\xi_k)
\simeq\bigg(\sum_{1\le s\le k}\Vert\xi_s\Vert^{2p/s}_h\bigg)^{1/p}
$$
(by Lemma 2.12, the errors can be taken arbitrary small in 
$C^\infty$ topology). We write $\xi_s=\sum_{1\le\alpha\le r}\xi_{s\alpha}
e_\alpha$. By (2.7) we have
$$
\Vert \xi_s\Vert_h^2=
\sum_\alpha|\xi_{s\alpha}|^2+
\sum_{i,j,\alpha,\beta}c_{ij\alpha\beta}z_i\overline z_j\xi_{s\alpha}
\overline \xi_{s\beta}+O(|z|^3|\xi|^2).
$$
The question is to evaluate the curvature of the weighted metric defined by
$$
\eqalign{
\Psi(z\,;\,\xi_1,\ldots,\xi_k)
&=\bigg(\sum_{1\le s\le k}\Vert\xi_s\Vert^{2p/s}_h\bigg)^{1/p}\cr
&=\bigg(\sum_{1\le s\le k}\Big(\sum_\alpha|\xi_{s\alpha}|^2+\!\!
\sum_{i,j,\alpha,\beta}c_{ij\alpha\beta}z_i\overline z_j\xi_{s\alpha}
\overline\xi_{s\beta}
\Big)^{p/s}\bigg)^{1/p}\kern-4pt{}+O(|z|^3).\cr}
$$
We set $|\xi_s|^2=\sum_\alpha|\xi_{s\alpha}|^2$. A straightforward 
calculation yields
$$
\eqalign{
&\log\Psi(z\,;\,\xi_1,\ldots,\xi_k)=\cr
&~~{}={1\over p}\log\sum_{1\le s\le k}|\xi_s|^{2p/s}+
\sum_{1\le s\le k}{1\over s}{|\xi_s|^{2p/s}\over \sum_t|\xi_t|^{2p/t}}
\sum_{i,j,\alpha,\beta}c_{ij\alpha\beta}z_i\overline z_j
{\xi_{s\alpha}\overline\xi_{s\beta}\over|\xi_s|^2}+O(|z|^3).\cr}
$$
By (2.11), the curvature form of $L_k=\cO_{X_k^\GG}(1)$ 
is given at the central point $x_0$ by the following formula.

\claim (2.13) Proposition|With the above choice of coordinates and in
terms of the rescaled components $\xi_s=\varepsilon_s^s\nabla^sf(0)$ at 
$x_0\in X$, we have the approximate expression
$$
\Theta_{L_k,\Psi^*_{h,p,\varepsilon}}(x_0,[\xi])\simeq
\omega_{a,r,p}(\xi)+{i\over 2\pi}
\sum_{1\le s\le k}{1\over s}{|\xi_s|^{2p/s}\over \sum_t|\xi_t|^{2p/t}}
\sum_{i,j,\alpha,\beta}c_{ij\alpha\beta}
{\xi_{s\alpha}\overline\xi_{s\beta}\over|\xi_s|^2}\,dz_i\wedge d\overline z_j
$$
where the error terms  are
$O(\max_{2\le s\le k}(\varepsilon_s/\varepsilon_{s-1})^s)$ uniformly on
the compact variety $X_k^\GG$. Here $\omega_{a,r,p}$ is the $($degenerate$)$
K\"ahler metric associated with the weight $a=(1^{[r]},2^{[r]},\ldots,k^{[r]})$ 
of the canonical $\bC^*$ action on $J^kV$.
\endclaim

Thanks to the uniform approximation, we can (and will) neglect the error 
terms in the calculations below. Since $\omega_{a,r,p}$ is positive definite
on the fibers of $X_k^\GG\to X$ (at least outside of the axes $\xi_s=0$), 
the index of the $(1,1)$
curvature form $\Theta_{L_k,\Psi^*_{h,p,\varepsilon}}(z,[\xi])$ is equal
to the index of the $(1,1)$-form
$$
\gamma_k(z,\xi):={i\over 2\pi}
\sum_{1\le s\le k}{1\over s}{|\xi_s|^{2p/s}\over \sum_t|\xi_t|^{2p/t}}
\sum_{i,j,\alpha,\beta}c_{ij\alpha\beta}(z)
{\xi_{s\alpha}\overline\xi_{s\beta}\over|\xi_s|^2}\,dz_i\wedge d\overline z_j
\leqno(2.14)
$$
depending only on the differentials $(dz_j)_{1\le j\le n}$ on~$X$. The 
$q$-index integral of $(L_k,\Psi^*_{h,p,\varepsilon})$ on $X^\GG_k$ is 
therefore equal to
$$
\eqalign{
&\int_{X^\GG_k(L_k,q)}\Theta_{L_k,\Psi^*_{h,p,\varepsilon}}^{n+kr-1}=\cr
&\qquad{}={(n+kr-1)!\over n!(kr-1)!}
\int_{z\in X}\int_{\xi\in P(1^{[r]},\ldots,k^{[r]})}
\omega_{a,r,p}^{kr-1}(\xi)\bone_{\gamma_k,q}(z,\xi)\gamma_k(z,\xi)^n\cr}
$$
where $\bone_{\gamma_k,q}(z,\xi)$ is the characteristic function of the open
set of points where $\gamma_k(z,\xi)$ has signature $(n-q,q)$ in terms of
the $dz_j$'s. Notice that since $\gamma_k(z,\xi)^n$ is~a determinant, the
product $\bone_{\gamma_k,q}(z,\xi)\gamma_k(z,\xi)^n$ gives rise to a continuous
function on~$X^\GG_k$. Formula 1.13 with $r_1=\ldots=r_k=r$ and
$a_s=s$ yields the slightly more explicit
integral
$$
\eqalign{
&\int_{X^\GG_k(L_k,q)}\Theta_{L_k,\Psi^*_{h,p,\varepsilon}}^{n+kr-1}=
{(n+kr-1)!\over n!(k!)^r}~~\times\cr
&\qquad\int_{z\in X}\int_{(x,u)\in\Delta_{k-1}\times(S^{2r-1})^k}
\bone_{g_k,q}(z,x,u)g_k(z,x,u)^n\,
{(x_1\ldots x_k)^{r-1}\over (r-1)!^k}\,dx\,d\mu(u),\cr}
$$
where $g_k(z,x,u)=\gamma_k(z,x_1^{1/2p}u_1,\ldots,x_k^{k/2p}u_k)$ is given by
$$
g_k(z,x,u)={i\over 2\pi}\sum_{1\le s\le k}{1\over s}x_s
\sum_{i,j,\alpha,\beta}c_{ij\alpha\beta}(z)\,
u_{s\alpha}\overline u_{s\beta}\,dz_i\wedge d\overline z_j
\leqno(2.15)
$$
and $\bone_{g_k,q}(z,x,u)$ is the characteristic function of its $q$-index 
set. Here 
$$
d\nu_{k,r}(x)=(kr-1)!\,{(x_1\ldots x_k)^{r-1}\over (r-1)!^k}\,dx
\leqno(2.16)
$$
is a probability measure on $\Delta_{k-1}$, and we can rewrite
$$
\leqalignno{
&\int_{X^\GG_k(L_k,q)}\Theta_{L_k,\Psi^*_{h,p,\varepsilon}}^{n+kr-1}
={(n+kr-1)!\over n!(k!)^r(kr-1)!}~~\times\cr
&\qquad\int_{z\in X}
\int_{(x,u)\in\Delta_{k-1}\times(S^{2r-1})^k}
\bone_{g_k,q}(z,x,u)g_k(z,x,u)^n\,d\nu_{k,r}(x)\,d\mu(u).&(2.17)\cr}
$$
Now, formula (2.15) shows that $g_k(z,x,u)$ is a ``Monte Carlo''
evaluation of the curvature tensor, obtained by averaging the curvature
at random points $u_s\in S^{2r-1}$ with certain positive weights $x_s/s\,$; 
we should then think of the \hbox{$k$-jet}
$f$ as some sort of random parameter such that the derivatives 
$\nabla^kf(0)$ are uniformly distributed in all directions. Let us compute
the expected value of
$(x,u)\mapsto g_k(z,x,u)$ with respect to the probability measure
$d\nu_{k,r}(x)\,d\mu(u)$. Since 
$\int_{S^{2r-1}}u_{s\alpha}\overline u_{s\beta}d\mu(u_s)={1\over r}
\delta_{\alpha\beta}$ and $\int_{\Delta_{k-1}}x_s\,d\nu_{k,r}(x)={1\over k}$,
we find
$$
{\bf E}(g_k(z,\bu,\bu))={1\over kr}
\sum_{1\le s\le k}{1\over s}\cdot{i\over 2\pi}\sum_{i,j,\alpha}
c_{ij\alpha\alpha}(z)\,dz_i\wedge d\overline z_j.
$$
In other words, we get the normalized trace of the curvature, i.e.
$$
{\bf E}(g_k(z,\bu,\bu))={1\over kr}
\Big(1+{1\over 2}+\ldots+{1\over k}\Big)\Theta_{\det(V^*),\det h^*},
\leqno(2.18)
$$
where $\Theta_{\det(V^*),\det h^*}$ is the $(1,1)$-curvature form of
$\det(V^*)$ with the metric induced by~$h$. It is natural to guess that 
$g_k(z,x,u)$ behaves asymptotically as its expected value
${\bf E}(g_k(z,\bu,\bu))$ when $k$ tends to infinity. If we replace brutally 
$g_k$ by its expected value in (2.17), we get the integral
$$
{(n+kr-1)!\over n!(k!)^r(kr-1)!}{1\over (kr)^n}
\Big(1+{1\over 2}+\ldots+{1\over k}\Big)^n\int_X\bone_{\eta,q}\eta^n,
$$
where $\eta:=\Theta_{\det(V^*),\det h^*}$ and $\bone_{\eta,q}$ is the
characteristic function of its $q$-index set in~$X$. The leading constant is
equivalent to $(\log k)^n/n!(k!)^r$ modulo 
a multiplicative factor $1+O(1/\log k)$. By working out a more precise analysis
of the deviation, we will prove the following result.

\claim (2.19) Probabilistic estimate|Fix smooth hermitian metrics $h$ on $V$ and
$\omega={i\over 2\pi} \sum\omega_{ij}dz_i\wedge d\overline z_j$ on $X$. 
Denote by $\Theta_{V,h}=-{i\over 2\pi}\sum
c_{ij\alpha\beta}dz_i\wedge d\overline z_j\otimes e_\alpha^*\otimes
e_\beta$ the curvature tensor of $V$ with respect to an $h$-orthonormal frame
$(e_\alpha)$, and put
$$
\eta(z)=\Theta_{\det(V^*),\det h^*}={i\over 2\pi}\sum_{1\le i,j\le n}\eta_{ij}
dz_i\wedge d\overline z_j,\qquad
\eta_{ij}=\sum_{1\le\alpha\le r}c_{ij\alpha\alpha}.
$$
Finally consider the $k$-jet line bundle $L_k=\smash{\cO_{X_k^\GG}(1)}\to
X_k^\GG$ equipped with the induced metric $\Psi^*_{h,p,\varepsilon}$
$($as defined above, with $1=\varepsilon_1\gg\varepsilon_2\gg\ldots\gg
\varepsilon_k>0)$. When $k$ tends 
to infinity, the integral of the top power of the curvature of $L_k$ on its
$q$-index set $X^\GG_k(L_k,q)$ is given by
$$
\int_{X^\GG_k(L_k,q)}\Theta_{L_k,\Psi^*_{h,p,\varepsilon}}^{n+kr-1}=
{(\log k)^n\over n!\,(k!)^r}\bigg(
\int_X\bone_{\eta,q}\eta^n+O((\log k)^{-1})\bigg)
$$
for all $q=0,1,\ldots,n$, and the error term $O((\log k)^{-1})$ can be 
bounded explicitly in terms of $\Theta_V$, $\eta$ and $\omega$. Moreover, the 
left hand side is identically zero for $q>n$.
\endclaim

The final statement follows from the observation that the curvature of
$L_k$ is positive along the fibers of $X_k^\GG\to X$, by the 
plurisubharmonicity of the weight (this is true even 
when the partition of unity terms are taken into account, since they
depend only on the base); therefore the $q$-index sets are empty for
$q>n$. We start with three elementary lemmas.

\claim (2.20) Lemma| The integral
$$I_{k,r,n}=\int_{\Delta_{k-1}}
\bigg(\sum_{1\le s\le k}{x_s\over s}\bigg)^nd\nu_{k,r}(x)$$
is given by the expansion
$$
I_{k,r,n}=\sum_{1\le s_1,s_2,\ldots,s_n\le k}{1\over s_1s_2\ldots s_n}\;
{(kr-1)!\over  (r-1)!^k}\;{\prod_{1\le i\le k} (r-1+\beta_i)!\over(kr+n-1)!}.
\leqno{\rm(a)}
$$
where $\beta_i=\beta_i(s)=\card\{j\,;\;s_j=i\}$, $\sum\beta_i=n$, 
$1\le i\le k$. The quotient
$$
I_{k,r,n}\bigg/{r^n\over kr(kr+1)\ldots (kr+n-1)}
\Big(1+{1\over 2}+\ldots+{1\over k}\Big)^n
$$
is bounded below by $1$ and bounded above by
$$
1+{1\over 3}\sum_{m=2}^n{2^mn!\over (n-m)!}
\bigg(1+{1\over 2}+\ldots+{1\over k}\bigg)^{-m}=1+O((\log k)^{-2})
\leqno{\rm(b)}
$$
As a consequence
$$
\leqalignno{
I_{k,r,n}&={1\over k^n}\Big(\Big(1+{1\over 2}+\ldots+{1\over k}\Big)^n+
O((\log k)^{n-2})\Big)&{\rm(c)}\cr
&={(\log k+\gamma)^n+O((\log k)^{n-2})\over k^n}
\cr}
$$
where $\gamma$ is the Euler-Mascheroni constant.
\endclaim

\proof. Let us expand the $n$-th power 
$\big(\sum_{1\le s\le k}{x_s\over s}\big)^n$. This gives
$$
I_{k,r,n}=\sum_{1\le s_1,s_2,\ldots,s_n\le k}{1\over s_1s_2\ldots s_n}
\int_{\Delta_{k-1}}
x_1^{\beta_1}\ldots x_k^{\beta_k}\,d\nu_{k,r}(x)$$
and by definition of the measure $\nu_{k,r}$ we have
$$\int_{\Delta_{k-1}}
x_1^{\beta_1}\ldots x_k^{\beta_k}\,d\nu_{k,r}(x)=
{(kr-1)!\over (r-1)!^k}\int_{\Delta_{k-1}}
x_1^{r+\beta_1-1}\ldots x_k^{r+\beta_k-1}\,dx_1\ldots dx_k.
$$
By Formula (1.14), we find
$$\eqalign{
\int_{\Delta_{k-1}}
x_1^{\beta_1}\ldots x_k^{\beta_k}\,d\nu_{k,r}(x)&=
{(kr-1)!\over (r-1)!^k}\;{\prod_{1\le i\le k}(r+\beta_i-1)!\over (kr+n-1)!}\cr
&={r^n\prod_{i,\,\beta_i\ge 1}(1+{1\over r})(1+{2\over r})\ldots
(1+{\beta_i-1\over r})\over kr(kr+1)\ldots(kr+n-1)},\cr}
$$
and (2.20$\,$a) follows from the first equality. The final product is 
minimal \hbox{when $r=1$}, thus
$$
\leqalignno{
&{r^n\over kr(kr+1)\ldots(kr+n-1)}\le
\int_{\Delta_{k-1}}
x_1^{\beta_1}\ldots x_k^{\beta_k}\,d\nu_{k,r}(x)\cr
&\kern7cm\le{r^n\prod_{1\le i\le k}\beta_i!\over kr(kr+1)\ldots(kr+n-1)}.
&(2.21)\cr}
$$
Also, the integral is maximal when all $\beta_i$ vanish except one, in 
which case one gets
$$\int_{\Delta_{k-1}}
x_j^n\,d\nu_{k,r}(x)=
{r(r+1)\ldots(r+n-1)\over kr(kr+1)\ldots (kr+n-1)}.
\leqno(2.22)
$$
By (2.21), we find the lower and upper bounds
$$
\leqalignno{
I_{k,r,n}&\ge {r^n\over kr(kr+1)\ldots (kr+n-1)}
\Big(1+{1\over 2}+\ldots+{1\over k}\Big)^n,&(2.23)\cr
I_{k,r,n}&\le {r^n\over kr(kr+1)\ldots (kr+n-1)}
\sum_{1\le s_1,\ldots,s_n\le k}{\beta_1!\ldots \beta_k!\over s_1\ldots s_n}.
&(2.24)\cr}
$$
In order to make the upper bound more explicit, we reorganize the
$n$-tuple $(s_1,\ldots,s_n)$ into those indices $t_1<\ldots<t_\ell$
which appear a certain number of times $\alpha_i=\beta_{t_i}\ge 2$, and 
those, say
$t_{\ell+1}<\ldots<t_{\ell+m}$, which appear only once. We have of
course $\sum\beta_i=n-m$, and each choice of the $t_i$'s corresponds
to $n!/\alpha_1!\ldots\alpha_\ell!$ possibilities for the $n$-tuple
$(s_1,\ldots,s_n)$. Therefore we get
$$
\sum_{1\le s_1,\ldots,s_n\le k}{\beta_1!\ldots \beta_k!\over s_1\ldots s_n}\le
n!\sum_{m=0}^n~~\sum_{\ell,\,\Sigma\alpha_i=n-m}~~\sum_{(t_i)}
{1\over t_1^{\alpha_1}\ldots t_\ell^{\alpha_\ell}}{1\over t_{\ell+1}
\ldots t_{\ell+m}}.
$$
A trivial comparison series vs.\ integral yields
$$
\sum_{s<t<+\infty}{1\over t^\alpha}\le {1\over \alpha-1}{1\over s^{\alpha-1}}
$$
and in this way, using successive integrations in $t_\ell$, $t_{\ell-1}$,
$\ldots$~, we get inductively
$$
\sum_{1\le t_1<\ldots<t_\ell<+\infty}{1\over t_1^{\alpha_1}\ldots t_\ell^{\alpha_\ell}}
\le{1\over \prod_{1\le i\le\ell}(\alpha_{\ell-i+1}+\ldots+\alpha_\ell-i)}
\le {1\over \ell!}\,,
$$
since $\alpha_i\ge 2$ implies $\alpha_{\ell-i+1}+\ldots+\alpha_\ell-i\ge i$. On 
the other hand
$$
\sum_{1\le t_{\ell+1}<\ldots<t_{\ell+m}\le k}{1\over t_{\ell+1}\ldots t_{\ell+m}}
\le{1\over m!}\kern-2pt\sum_{1\le s_1,\ldots,s_m\le k}\kern-2pt
{1\over s_1\ldots s_m}={1\over m!}\bigg(1+{1\over 2}+\ldots+{1\over k}\bigg)^m.
$$
Since partitions $\alpha_1+\ldots+\alpha_\ell=n-m$ satisfying the additional
restriction $\alpha_i\ge 2$ correspond to $\alpha_i'=\alpha_i-2$ satisfying 
$\sum \alpha'_i=n-m-2\ell$, their number is equal to
$$
{n-m-2\ell+\ell-1\choose \ell-1}={n-m-\ell-1\choose \ell-1}\le
2^{n-m-\ell-1}
$$
and we infer from this
$$
\sum_{1\le s_1,\ldots,s_n\le k}{\beta_1!\ldots \beta_k!\over s_1\ldots s_n}\le
\sum_{\ell\ge 1\atop 2\ell+m\le n}{2^{n-m-\ell-1}n!\over \ell!\,m!}
\bigg(1+{1\over 2}+\ldots+{1\over k}\bigg)^m
+\bigg(1+{1\over 2}+\ldots+{1\over k}\bigg)^n
$$
where the last term corresponds to the special case $\ell=0$, $m=n$. Therefore
$$
\eqalign{
\sum_{1\le s_i\le k}{\beta_1!\ldots \beta_k!\over s_1\ldots s_n}
&\le{e^{1/2}-1\over 2}\sum_{m=0}^{n-2}{2^{n-m}n!\over m!}
\bigg(1{+}{1\over 2}{+}\ldots{+}{1\over k}\bigg)^m\!+
\bigg(1{+}{1\over 2}{+}\ldots{+}{1\over k}\bigg)^n\cr
&\le {1\over 3}\sum_{m=2}^n{2^mn!\over (n-m)!}
\bigg(1{+}{1\over 2}{+}\ldots{+}{1\over k}\bigg)^{n-m}\!+
\bigg(1{+}{1\over 2}{+}\ldots{+}{1\over k}\bigg)^n.\cr}
$$
This estimate combined with (2.23,~2.24) implies the upper bound
(2.20~b) (the lower bound $1$ being now obvious). The asymptotic estimate
(2.20~c) follows immediately.\qed

\claim (2.25) Lemma|If $A$ is a hermitian $n\times n$ matrix, set
$\bone_{A,q}$ to be equal to $1$ if $A$ has signature $(n-q,q)$ and
$0$ otherwise. Then for all $n\times n$ hermitian matrices $A$,~$B$
we have the estimate
$$
\big|\bone_{A,q}\det A-\bone_{B,q}\det B\big|\le \Vert A-B\Vert~
\sum_{0\le i\le n-1}\Vert A\Vert^i\Vert B\Vert^{n-1-i},
$$
where $\Vert A\Vert$, $\Vert B\Vert$ are the hermitian operator norms 
of the matrices.
\endclaim

\proof. We first check that the estimate holds true for $|\det A-\det B|$.
Let $\lambda_1\le\ldots\le \lambda_n$ be the eigenvalues of $A$
and $\lambda'_1\le\ldots\le \lambda'_n$ be the eigenvalues of $B$.
We have $|\lambda_i|\le \Vert A\Vert$, $|\lambda'_i|\le \Vert B\Vert$
and the minimax principle implies that 
$|\lambda_i-\lambda'_i|\le\Vert A-B\Vert$.
We then get the desired estimate by writing
$$
\det A-\det B=\lambda_1\ldots\lambda_n-\lambda'_1\ldots\lambda'_n=
\sum_{1 \le i\le n}\lambda_1\ldots\lambda_{i-1}(\lambda_i-\lambda_i')
\lambda'_{i+1}\ldots\lambda'_n.
$$
This already implies (2.25) if $A$ or $B$ is degenerate. If $A$ and $B$ are 
non degenerate we only have to prove the result
when one of them (say $A$) has signature $(n-q,q)$ and the other
one (say $B$) has a different signature. If we put $M(t)=(1-t)A+tB$, the 
already established estimate for the determinant yields 
$$\Big|{d\over dt}\det M(t)\Big|\le n\Vert A-B\Vert\;\Vert M(t)\Vert
\le n\Vert A-B\Vert\big((1-t)\Vert A\Vert+t\Vert B\Vert\big)^{n-1}.
$$
However, since the signature of $M(t)$ is not the same for $t=0$ and $t=1$, 
there must exist $t_0\in{}]0,1[$ such that $(1-t_0)A+t_0B$ is degenerate.
Our claim follows by integrating the differential estimate on
the smallest such interval $[0,t_0]$, after observing that $M(0)=A$, 
$\det M(t_0)=0$, and that the integral of the right hand side on $[0,1]$ 
is the announced bound.\qed

\claim (2.26) Lemma|Let $Q_A$ be the hermitian quadratic form associated with
the hermitian operator $A$ on $\bC^n$. If $\mu$ is the rotation 
invariant probability measure on the unit sphere $S^{2n-1}$ of $\bC^n$ 
and $\lambda_i$ are the eigenvalues of $A$, we have
$$
\int_{|\zeta|=1}|Q_A(\zeta)|^2d\mu(\zeta)={1\over n(n+1)}\Big(
\sum \lambda_i^2+\Big(\sum\lambda_i\Big)^2\Big).
$$
The norm $\Vert A\Vert=\max|\lambda_i|$ satisfies the estimate
$$
{1\over n^2}\Vert A\Vert^2 
\le\int_{|\zeta|=1}|Q_A(\zeta)|^2d\mu(\zeta)\le\Vert A\Vert^2.
$$
\endclaim

\proof. The first identity as an easy calculation, and the inequalities 
follow by computing the eigenvalues of the quadratic form 
$\sum \lambda_i^2+\big(\sum\lambda_i\big)^2-c\lambda_{i_0}^2$, $c>0$.
The lower bound is attained e.g.\
for $Q_A(\zeta)=|\zeta_1|^2-{1\over n}(|\zeta_2|^2+\ldots+|\zeta_n|^2)$
when we take $i_0=1$ and $c=1+{1\over n}$.\qed

\proof\ {\rm of Proposition 2.19}. Take a vector 
$\zeta\in T_{X,z}$, $\zeta=\sum\zeta_i{\partial\over\partial z_i}$,
with~$\Vert\zeta\Vert_\omega=1$, and introduce the trace free sesquilinear
quadratic form
$$Q_{z,\zeta}(u)=\sum_{i,j,\alpha,\beta}
\wt c_{ij\alpha\beta}(z)\,\zeta_i\overline\zeta_j\,u_\alpha
\overline u_\beta,\qquad
\wt c_{ij\alpha\beta} = c_{ij\alpha\beta} -
{1\over r}\eta_{ij}\delta_{\alpha\beta},\qquad u\in\bC^r
$$
where $\eta_{ij}=\sum_{1\le\alpha\le r}c_{ij\alpha\alpha}$. We consider the corresponding 
trace free curvature tensor
$$
\wt\Theta_V={i\over 2\pi}\sum_{i,j,\alpha,\beta}
\wt c_{ij\alpha\beta}\,dz_i\wedge d\overline z_j\otimes e_\alpha^*\otimes e_\beta.
\leqno(2.27)
$$
As a general matter of notation, we adopt here the convention that the cano\-nical 
correspondence between hermitian forms and $(1,1)$-forms is 
normalized as $\sum a_{ij}dz_i\otimes d\overline z_j
\leftrightarrow {i\over 2\pi}\sum a_{ij}dz_i\wedge d\overline z_j$, and we take
the liberty of using the same symbols for both types of objects; we do so especially 
for $g_k(z,x,u)$ and $\eta(z)={i\over 2\pi}\sum \eta_{ij}(z)dz_i\wedge d\overline z_j=
\Tr \Theta_V(z)$. First observe that for all $k$-tuples of unit vectors 
$u=(u_1,\ldots,u_k)\in (S^{2r-1})^k$, $u_s=(u_{s\alpha})_{1\le\alpha\le r}$, we have
$$
\int_{(S^{2r-1})^k}\bigg|\sum_{1\le s\le k}{1\over s}x_s
\sum_{i,j,\alpha,\beta}\wt c_{ij\alpha\beta}(z)\,\zeta_i\overline\zeta_j
u_{s\alpha}\overline u_{s\beta}\bigg|^2d\mu(u)
=\sum_{1\le s\le k}{x_s^2\over s^2}{\bf V}(Q_{z,\zeta})
$$
where ${\bf V}(Q_{z,\zeta})$ is the variance of $Q_{z,\zeta}$
on $S^{2r-1}$. This is so because we have a sum over $s$ of
independent random variables on $(S^{2r-1})^k$, all of which have zero 
mean value. (Lemma 2.26 shows that the variance ${\bf V}(Q)$
of a  trace free hermitian quadratic form 
$Q(u)=\sum_{1\le\alpha\le r}\lambda_\alpha|u_\alpha|^2$ on the unit 
sphere $S^{2r-1}$ is equal to $\smash{{1\over r(r+1)}\sum\lambda_\alpha^2}\,$,
but we only give the formula to fix the ideas). Formula (2.22) yields
$$
\int_{\Delta_{k-1}}x_s^2d\nu_{k,r}(x)={r+1\over k(kr+1)}.
$$
Therefore, according to notation (2.15), we obtain 
the partial variance formula
$$
\eqalign{
\int_{\Delta_{k-1}\times (S^{2r-1})^k}\big|g_k(z,x,u)(\zeta)
&-\overline g_k(z,x)(\zeta)|^2d\nu_{k,r}(x)d\mu(u)\cr
&={(r+1)\over k(kr+1)}\bigg(\sum_{1\le s\le k}{1\over s^2}\bigg)
\sigma_h(\wt\Theta_V(\zeta,\zeta))^2\cr}
$$
in which
$$
\eqalign{
\overline g_k(z,x)(\zeta)&=\sum_{1\le s\le k}{1\over s}x_s{1\over r}
\sum_{ij\alpha}c_{ij\alpha\alpha}\zeta_i\overline\zeta_j=
\bigg(\sum_{1\le s\le k}{1\over s}x_s\bigg){1\over r}\eta(z)(\zeta),\cr
\sigma_h(\wt\Theta_V(\zeta,\zeta))^2&=
{\bf V}\big(u\mapsto\langle\wt\Theta_V(\zeta,\zeta)u,u\rangle_h\big)
=\int_{u\in S^{2r-1}}
\big|\langle\wt\Theta_V(\zeta,\zeta)u,u\rangle_h\big|^2d\mu(u).
\cr}
$$
By integrating over $\zeta\in S^{2n-1}\subset \bC^n$ and applying the 
left hand inequality in Lemma 2.26 we infer
$$
\leqalignno{
\int_{\Delta_{k-1}\times (S^{2r-1})^k}\big\Vert g_k(z,x,u)-
 \overline g_k(z,x)\Vert_\omega^2&d\nu_{k,r}(x)d\mu(u)\cr
&\le {n^2(r+1)\over k(kr+1)}\bigg(\sum_{1\le s\le k}{1\over s^2}\bigg)
\sigma_{\omega,h}(\wt\Theta_V)^2
&(2.28)\cr}
$$
where $\sigma_{\omega,h}(\wt\Theta_V)$ is the standard deviation of
$\wt\Theta_V$ on $S^{2n-1}\times S^{2r-1}$~:
$$
\sigma_{\omega,h}(\wt\Theta_V)^2=\int_{|\zeta|_\omega=1,\,
|u|_h=1}\big|\langle\wt\Theta_V(\zeta,\zeta)u,u\rangle_h\big|^2d\mu(\zeta)\,
d\mu(u).
$$
On the other hand, brutal estimates give the hermitian operator norm
estimates
$$
\leqalignno{
\Vert\overline g_k(z,x)\Vert_\omega&\le 
\bigg(\sum_{1\le s\le k}{1\over s}x_s\bigg)
{1\over r}\Vert\eta(z)\Vert_\omega,&(2.29)\cr
\Vert g_k(z,x,u)\Vert_\omega&\le \bigg(\sum_{1\le s\le k}{1\over s}x_s
\bigg)\Vert\Theta_V\Vert_{\omega,h}&(2.30)\cr}
$$
where
$$
\Vert\Theta_V\Vert_{\omega,h}=\sup_{|\zeta|_\omega=1,\,|u|_h=1}
\big|\langle\Theta_V(\zeta,\zeta)u,u\rangle_h\big|.
$$
We use these estimates to evaluate the $q$-index integrals. The integral 
associated with $\overline g_k(z,x)$ is much easier to 
deal with than $g_k(z,x,u)$ since the characteristic function of
the $q$-index set depends only on $z$. By Lemma 2.25 we find
$$
\eqalign{
&\big|\bone_{g_k,q}(z,x,u)\det g_k(z,x,u)-
\bone_{\eta,q}(z)\det\overline g_k(z,x)\big|\cr
&\qquad{}\le
\big\Vert g_k(z,x,u)-\overline g_k(z,x)\big\Vert_\omega\sum_{0\le i\le n-1}
\Vert g_k(z,x,u)\Vert_\omega^i\Vert\overline g_k(z,x)\Vert_\omega^{n-1-i}.\cr}
$$
The Cauchy-Schwarz inequality combined with (2.28 -- 2.30) implies
$$
\eqalign{
&\int_{\Delta_{k-1}\times (S^{2r-1})^k}
\big|\bone_{g_k,q}(z,x,u)\det g_k(z,x,u)-
\bone_{\eta,q}(z)\det\overline g_k(z,x)\big|\,d\nu_{k,r}(x)d\mu(u)\cr
&\quad{}\le\bigg(\int_{\Delta_{k-1}\times (S^{2r-1})^k}\big\Vert g_k(z,x,u)-
\overline g_k(z,x)\big\Vert_\omega^2d\nu_{k,r}(x)d\mu(u)\bigg)^{1/2}\times\cr
&\kern22pt\bigg(\int_{\Delta_{k-1}\times (S^{2r-1})^k}
\kern-2pt\bigg(\sum_{0\le i\le n-1}\kern-4pt
\Vert g_k(z,x,u)\Vert_\omega^i\Vert\overline g_k(z,x)
\Vert_\omega^{n-1-i}\bigg)^2d\nu_{k,r}(x)d\mu(u)\kern-2pt\bigg)^{1/2}\cr
&\quad{}\le {n(1+1/r)^{1/2}\over (k(k+1/r))^{1/2}}
\bigg(\sum_{1\le s\le k}{1\over s^2}\bigg)^{1/2}\sigma_{\omega,h}(\wt\Theta_V)
\sum_{1\le i\le n-1}\!\Vert\Theta_V\Vert_{\omega,h}^i
\Big({1\over r}\Vert\eta(z)\Vert_\omega\Big)^{n-1-i}\cr
&\kern45pt{}\times\bigg(\int_{\Delta_{k-1}}
\bigg(\sum_{1\le s\le k}{x_s\over s}\bigg)^{2n-2}d\nu_{k,r}(x)\bigg)^{1/2}
{}={}~~O\Big({(\log k)^{n-1}\over k^n}\Big)\cr}
$$
by Lemma 2.20 with $n$ replaced by $2n-2$. This is the essential error 
estimate. As one can see, the growth of the error 
mainly depends on the final integral factor, since the initial
multiplicative factor is
uniformly bounded over $X$. In order to get the principal term, we compute
$$
\eqalign{
\int_{\Delta_{k-1}}\det \overline g_k(z,x)\,d\nu_{k,r}(x)
&={1\over r^n}\det\eta(z)\int_{\Delta_{k-1}}
\bigg(\sum_{1\le s\le k}{x_s\over s}\bigg)^nd\nu_{k,r}(x)\cr
&\sim{(\log k)^n\over r^nk^n}\det\eta(z).\cr}
$$
From there we conclude that
$$
\eqalign{
\int_{z\in X}\int_{(x,u)\in\Delta_{k-1}\times (S^{2r-1})^k}
\bone_{g_k,q}(z,x,u) &g_k(z,x,u)^n\,d\nu_{k,r}(x)d\mu(u)\cr
&={(\log k)^n\over r^nk^n}
\int_X\bone_{\eta,q}\eta^n+O\Big({(\log k)^{n-1}\over k^n}\Big)\cr}
$$
The probabilistic estimate 2.19 follows by (2.17).\qed

\claim (2.31) Remark|{\rm If we take care of the precise bounds obtained
above, the proof gives in fact the explicit estimate
$$
\int_{X^\GG_k(L_k,q)}\Theta_{L_k,\Psi^*_{h,p,\varepsilon}}^{n+kr-1}=
{(n+kr-1)!\;I_{k,r,n}\over n!(k!)^r(kr-1)!}\bigg(
\int_X\bone_{\eta,q}\eta^n+\varepsilon_{k,r,n}J\bigg)
$$
where
$$
J=n\,(1+1/r)^{1/2}\bigg(\sum_{s=1}^k{1\over s^2}\bigg)^{1/2}
\int_X\sigma_{\omega,h}(\wt\Theta_V)\sum_{i=1}^{n-1}
r^{i+1}\Vert\Theta_V\Vert_{\omega,h}^i
\Vert\eta(z)\Vert_\omega^{n-1-i}\omega^n
$$
and
$$
\eqalign{
|\varepsilon_{k,r,n}|&\le
{\displaystyle\bigg(\int_{\Delta_{k-1}}
\bigg(\sum_{s=1}^k{x_s\over s}\bigg)^{2n-2}d\nu_{k,r}(x)\bigg)^{1/2}\over
\displaystyle (k(k+1/r))^{1/2}\int_{\Delta_{k-1}}\phantom{\Bigg|}
\bigg(\sum_{s=1}^k{x_s\over s}\bigg)^nd\nu_{k,r}(x)}\cr
&\le
{\Big(1+{1\over 3}\sum_{m=2}^{2n-2}{2^m(2n-2)!\over(2n-2-m)!}
\big(1+{1\over 2}+\ldots+{1\over k}\big)^{-m}\Big)^{1/2}
\over 1+{1\over 2}+\ldots+{1\over k}}\sim{1\over \log k}\cr}
$$
by the lower and upper bounds of $I_{k,r,n}$, $I_{k,r,2n-2}$ obtained 
in Lemma 2.20. As $(2n-2)!/(2n-2-m)!\le (2n-2)^m$, one easily shows that
$$
|\varepsilon_{k,r,n}|\le {(31/15)^{1/2}\over \log k}\qquad\hbox{for $k\ge e^{5n-5}$}.\leqno(2.32)
$$
Also, we see that the error terms vanish if $\smash{\wt\Theta_V}$ is 
identically zero, but this is of course a rather unexpected 
circumstance. In general, since the form $\wt\Theta_V$ is trace free, 
Lemma~2.23 applied to the quadratic form
$u\mapsto\langle\wt\Theta_V(\zeta,\zeta)u,u\rangle$ on $\bC^r$ implies 
$\sigma_{\omega,h}(\wt\Theta_V)\le (r+1)^{-1/2}\Vert\wt\Theta_V\Vert_{\omega,h}$.
This yields the simpler bound
$$\kern15pt
J\le n\,r^{1/2}\bigg(\sum_{s=1}^k{1\over s^2}\bigg)^{1/2}
\int_X\Vert\wt\Theta_V\Vert_{\omega,h}\sum_{i=1}^{n-1}
r^i\Vert\Theta_V\Vert_{\omega,h}^i
\Vert\eta(z)\Vert_\omega^{n-1-i}\omega^n.
\leqno\hbox to 0.1pt{\rlap{\rlap{(2.33)}\kern\hsize\llap{\square}}}
$$}
\endclaim

It will be useful to extend the above estimates to the case of sections of
$$
L_k=\cO_{X_k^\GG}(1)\otimes
\pi_k^*\cO\Big({1\over kr}\Big(1+{1\over 2}+\ldots+{1\over k}\Big)F\Big)
\leqno(2.34)
$$
where $F\in\Pic_\bQ(X)$ is an arbitrary $\bQ$-line bundle on~$X$ and 
$\pi_k:X_k^\GG\to X$ is the natural projection. We assume here
that $F$ is also equipped with a smooth hermitian metric $h_F$. In formula
(2.20), the renormalized metric $\eta_k(z,x,u)$ of $L_k$ takes the form
$$
\eta_k(z,x,u)={1\over {1\over kr}(1+{1\over 2}+\ldots+{1\over k})}g_k(z,x,u)+
\Theta_{F,h_F}(z),
\leqno(2.35)
$$
and by the same calculations its  expected value is
$$
\eta(z):={\bf E}(\eta_k(z,\bu,\bu))=\Theta_{\det V^*,\det h^*}(z)+
\Theta_{F,h_F}(z).
\leqno(2.36)
$$
Then the variance estimate for $\eta_k-\eta$ is unchanged, and the
$L^p$ bounds for $\eta_k$ are still valid, since our forms are just shifted
by adding the constant smooth term $\Theta_{F,h_F}(z)$. The probabilistic
estimate 2.18 is therefore still true in exactly the same form, provided
we use (2.34 -- 2.36) instead of the previously defined $L_k$, $\eta_k$
and~$\eta$. An application of holomorphic Morse inequalities gives the 
desired cohomology estimates for 
$$
\eqalign{
h^q\Big(X,E_{k,m}^\GG V^*&{}\otimes
\cO\Big({m\over kr}\Big(1+{1\over 2}+\ldots+{1\over k}\Big)F\Big)\Big)\cr
&{}=h^q(X_k^\GG,\cO_{X_k^\GG}(m)\otimes
\pi_k^*\cO\Big({m\over kr}\Big(1+{1\over 2}+\ldots+{1\over k}\Big)F\Big)\Big),
\cr}
$$
provided $m$ is sufficiently divisible to give a multiple of $F$ which
is a $\bZ$-line bundle.

\claim (2.37) Theorem|Let $(X,V)$ be a directed manifold, $F\to X$ a
$\bQ$-line bundle, $(V,h)$ and $(F,h_F)$ smooth hermitian structure on $V$ 
and $F$ respectively. We define
$$
\eqalign{
L_k&=\cO_{X_k^\GG}(1)\otimes
\pi_k^*\cO\Big({1\over kr}\Big(1+{1\over 2}+\ldots+{1\over k}\Big)F\Big),\cr
\eta&=\Theta_{\det V^*,\det h^*}+\Theta_{F,h_F}.\cr}
$$
Then for all $q\ge 0$ and all $m\gg k\gg 1$ such that 
m is sufficiently divisible, we have
$$\leqalignno{\kern20pt
h^q(X_k^\GG,\cO(L_k^{\otimes m}))&\le {m^{n+kr-1}\over (n+kr-1)!}
{(\log k)^n\over n!\,(k!)^r}\bigg(
\int_{X(\eta,q)}(-1)^q\eta^n+O((\log k)^{-1})\bigg),&\hbox{\rm(a)}\cr
h^0(X_k^\GG,\cO(L_k^{\otimes m}))&\ge {m^{n+kr-1}\over (n+kr-1)!}
{(\log k)^n\over n!\,(k!)^r}\bigg(
\int_{X(\eta,\le 1)}\eta^n-O((\log k)^{-1})\bigg),&\hbox{\rm(b)}\cr
\cr
\chi(X_k^\GG,\cO(L_k^{\otimes m}))&={m^{n+kr-1}\over (n+kr-1)!}
{(\log k)^n\over n!\,(k!)^r}\big(
c_1(V^*\otimes F)^n+O((\log k)^{-1})\big).&\hbox{\rm(c)}\cr
\cr}
$$
\vskip-4pt
\endclaim

Green and Griffiths [GG79] already checked the Riemann-Roch
calculation (2.37$\,$c) in the special case
$V=T^*_X$ and $F=\cO_X$. Their proof is much simpler since it relies only
on Chern class calculations, but it cannot provide any information on
the individual cohomology groups, except in very special cases where
vanishing theorems can be applied; in fact in dimension 2, the
Euler characteristic satisfies $\chi=h^0-h^1+h^2\le h^0+h^2$, hence
it is enough to get the vanishing of the top cohomology group $H^2$
to infer $h^0\ge\chi\,$; this works for surfaces by means of a well-known
vanishing theorem of Bogomolov which implies in general
$$H^n\bigg(X,E_{k,m}^\GG T_X^*\otimes
\cO\Big({m\over kr}\Big(1+{1\over 2}+\ldots+{1\over k}\Big)F\Big)\Big)\bigg)=0
$$
as soon as $K_X\otimes F$ is big and $m\gg 1$.

In fact, thanks to Bonavero's singular holomorphic Morse inequalities 
[Bon93], everything works almost unchanged 
in the case where $V\subset T_X$
has singularities and $h$ is an admissible metric on $V$ (see (0.11)).
We only have to find a blow-up $\mu:\smash{\wt X}_k\to X_k$ so that
the resulting pull-backs $\mu^*L_k$ and $\mu^*V$ are locally free,
and $\mu^*\det h^*$, $\mu^*\Psi_{h,p,\varepsilon}$ only have divisorial
singularities. Then $\eta$ is a $(1,1)$-current with logarithmic poles,
and we have to deal with smooth metrics on $
\mu^*L_k^{\otimes m}\otimes\cO(-mE_k)$ where $E_k$ is a certain effective 
divisor on $X_k$ (which, by our assumption (0.11), does not project onto
$X$). The cohomology groups involved are then the twisted
cohomology groups
$$
H^q(X_k^\GG,\cO(L_k^{\otimes m})\otimes\cJ_{k,m})
$$
where $\cJ_{k,m}=\mu_*(\cO(-mE_k))$ is the corresponding multiplier ideal sheaf,
and the Morse integrals need only be evaluated in the complement of the 
poles, that is on $X(\eta,q)\ssm S$ where $S=\Sing(V)\cup\Sing(h)$. Since
$$
(\pi_k)_*\big(\cO(L_k^{\otimes m})\otimes\cJ_{k,m}\big)\subset
E_{k,m}^\GG V^*\otimes
\cO\Big({m\over kr}\Big(1+{1\over 2}+\ldots+{1\over k}\Big)F\Big)\Big)
$$
we still get a lower bound for the $H^0$ of the latter sheaf (or for the $H^0$
of the un-twisted line bundle $\cO(L_k^{\otimes m})$ on $\smash{X_k^\GG}$).
If we assume that $K_V\otimes F$ is big, these considerations
also allow us to obtain a strong estimate in terms of the volume, by
using an approximate Zariski decomposition on a suitable blow-up of~$(X,V)$.
The following corollary implies in particular Theorem 0.5.

\claim (2.38) Corollary|If $F$ is an arbitrary $\bQ$-line bundle over~$X$, 
one has
$$
\eqalign{
h^0\bigg(&X_k^\GG,\cO_{X_k^\GG}(m)\otimes\pi_k^*\cO
\Big({m\over kr}\Big(1+{1\over 2}+\ldots+{1\over k}\Big)F\Big)\bigg)\cr
&\ge {m^{n+kr-1}\over (n+kr-1)!}
{(\log k)^n\over n!\,(k!)^r}\Big(
\Vol(K_V\otimes F)-O((\log k)^{-1})\Big)-o(m^{n+kr-1}),\cr}
$$
when $m\gg k\gg 1$, in particular there are many sections of the
$k$-jet differentials of degree $m$ twisted by the appropriate
power of $F$ if $K_V\otimes F$ is big.
\endclaim

\proof. The volume is computed here as usual, i.e.\ after performing a
suitable modifi\-cation $\mu:\smash{\wt X}\to X$ which converts $K_V$ into 
an invertible sheaf. There is of course nothing to prove if
$K_V\otimes F$ is not big, so we can assume $\Vol(K_V\otimes F)>0$.
Let us fix smooth hermitian metrics $h_0$ on $T_X$
and $h_F$ on $F$. They induce a metric $\mu^*(\det h_0^{-1}\otimes
h_F)$ on $\mu^*(K_V\otimes F)$ which, by our definition of $K_V$, is
a smooth metric. By the result of Fujita [Fuj94] on
approximate Zariski decomposition, for every $\delta>0$, one can find
a modification $\mu_\delta:\smash{\wt X_\delta}\to X$ dominating
$\mu$ such that
$$
\mu_\delta^*(K_V\otimes F) =\cO_{\wt X_\delta}(A+E)
$$
where $A$ and $E$ are $\bQ$-divisors, $A$ ample and $E$ effective,
with 
$$\Vol(A)=A^n\ge \Vol(K_V\otimes F)-\delta.$$
If we take a smooth metric $h_A$ with positive definite curvature form
$\Theta_{A,h_A}$, then we get a singular hermitian metric $h_Ah_E$ on
$\mu_\delta^*(K_V\otimes F)$ with poles along $E$, i.e.\ the quotient
$h_Ah_E/\mu^*(\det h_0^{-1}\otimes h_F)$ is of the form $e^{-\varphi}$ where
$\varphi$ is quasi-psh with log poles $\log|\sigma_E|^2$ 
(mod $C^\infty(\smash{\wt X_\delta}))$ precisely given
by the divisor~$E$. We then only need to take the singular metric $h$
on $T_X$ defined by
$$
h=h_0e^{{1\over r}(\mu_\delta)*\varphi}
$$
(the choice of the factor ${1\over r}$ is there to correct adequately 
the metric on $\det V$). By construction $h$ induces an 
admissible metric on $V$ and the resulting 
curvature current $\eta=\Theta_{K_V,\det h^*}+\Theta_{F,h_F}$ is such that
$$
\mu_\delta^*\eta = \Theta_{A,h_A} +[E],\qquad
\hbox{$[E]={}$current of integration on $E$.}
$$
Then the $0$-index Morse integral in the complement of the poles 
is given by
$$
\int_{X(\eta,0)\ssm S}\eta^n=\int_{\wt X_\delta}\Theta_{A,h_A}^n=A^n\ge
\Vol(K_V\otimes F)-\delta
$$
and (2.38) follows from the fact that $\delta$ can be taken arbitrary 
small.\qed

\claim (2.39) Example|{\rm In some simple cases, the above estimates can 
lead to very explicit results. Take for instance $X$ to be a smooth
complete intersection of multidegree $(d_1,d_2,\ldots,d_s)$ in $\bP^{n+s}_\bC$
and consider the absolute case $V=T_X$. Then 
$$K_X=\cO_X(d_1+\ldots+d_s-n-s-1).$$
Assume that $X$ is of general type, i.e.\ $\sum d_j>n+s+1$.
Let us equip $V=T_X$ with the restriction of the Fubini-Study metric 
$h=\Theta_{\cO(1)}\,$;
a better choice might be the K\"ahler-Einstein metric but we want to 
keep the calculations as elementary as possible. The standard formula for
the curvature tensor of a submanifold gives
$$
\Theta_{T_X,h}=(\Theta_{T_{\bP^{n+s}},h})_{|X}+\beta^*\wedge\beta
$$
where
$\smash{\beta\in 
C^\infty\big(\Lambda^{1,0}T^*_X\otimes\Hom(T_X,\bigoplus\cO(d_j))\big)}$
is the second fundamental form. In other words, by the well known
formula for the curvature of projective space, we have
$$
\langle\Theta_{T_X,h}(\zeta,\zeta)u,u\rangle=|\zeta|^2|u|^2
+|\langle\zeta,u\rangle|^2-|\beta(\zeta)\cdot u|^2.
$$
The curvature $\rho$ of $(K_X,\det h^*)$ (i.e.\ the opposite of the 
Ricci form $\Tr\Theta_{T_X,h}$) is given by
$$
\rho=-\Tr\Theta_{TX,h}=\Tr(\beta\wedge\beta^*)-(n+1)h\ge -(n+1)h.
\leqno(2.40)
$$
We take here $F=\cO_X(-a)$, $a\in\bQ_+$, and we want to determine 
conditions for the existence of sections
$$
H^0\bigg(X,E^\GG_{k,m}T^*_X\otimes\cO
\Big(-a{m\over kr}\Big(1+{1\over 2}+\ldots+{1\over k}\Big)\Big)\bigg),
\qquad m\gg 1.\leqno(2.41)
$$
We have to choose $K_X\otimes\cO_X(-a)$ ample, i.e.\
$\sum d_j>n+s+a+1$, and then (by an appropriate choice of the metric
of $F=\cO_X(-a)$), the form $\eta=\Theta_{K_X\otimes\cO_X(-a)}$ can be taken to
be any positive form cohomologous to $(\sum d_j-(n+s+a+1))h$.
We use remark 2.31 and estimate the error terms by considering
the K\"ahler metric
$$
\omega=\rho+(n+s+2)h\equiv \Big(\sum d_j+1\Big)h.
$$
Inequality (2.40) shows that $\omega\ge 2h$ and also that
$\omega\ge\Tr(\beta\wedge\beta^*)$. From this, one easily concludes 
that $\Vert\eta\Vert_\omega\le 1$ by an appropriate choice of $\eta$,
as well as $\Vert\Theta_{T_X,h}\Vert_{\omega,h}\le 1$ 
and $\Vert\wt\Theta_{T_X,h}\Vert_{\omega,h}\le 2$. By (2.33), we obtain
for $n\ge 2$
$$
J\le n^{3/2}{\pi\over \sqrt{6}}\times 2\,{n^n-1\over n-1}\int_X\omega^n
<{4\pi\over\sqrt{6}}\,n^{n+1/2}\int_X\omega^n
$$
where $\int_X\omega^n=\big(\sum d_j+1\big)^n\deg(X)$. On the other hand,
the leading term $\int_X\eta^n$ equals 
$\big(\sum d_j-n-s-a-1\big)^n\deg(X)$ with
$\deg(X)=d_1\ldots d_s$. By the bound (2.32 ) on the error term 
$\varepsilon_{k,r,n}$, we find that the leading coefficient of the
growth of our spaces of sections is strictly
controlled below by a multiple of
$$
\Big(\sum d_j-n-s-a-1\Big)^n-4\pi\Big({31\over 90}\Big)^{1/2}\,
{n^{n+1/2}\over \log k}\Big(\sum d_j+1\Big)^n
$$
if $k\ge e^{5n-5}$. A sufficient condition for the existence of
sections in (2.41) is thus
$$
k\ge\exp\Big(7.38\,n^{n+1/2}\Big({\sum d_j+1\over\sum d_j-n-s-a-1}\Big)^n\Big).
\leqno(2.42)
$$
This is good in view of the fact that we can cover arbitrary smooth 
complete intersections of general type. On the other hand, even when the
degrees $d_j$ tend to $+\infty$, we still get a large lower bound
$k\sim \exp(7.38\,n^{n+1/2})$ on the order of jets, and this is far 
from being optimal$\,$: Diverio [Div09] has shown e.g.\ that one can take
$k=n$ for smooth hypersurfaces of high degree. It is however not
unlikely that one could improve estimate (2.42) with more careful choices
of $\omega$, $h$.\qed}
\endclaim

\section{3. On the base locus of sections of 
{\twelvebfit k}-jet bundles}
\sectionrunning{3. On the base locus of sections of $k$-jet bundles}

The final step required for a complete solution of the Green-Griffiths
conjecture would be to calculate the base locus $B_k\subset X_k^{\GG}$
of the space of sections
$$
H^0(X_k^\GG,\cO_{X_k^\GG}(m)\otimes\pi_k^*\cO(-m\delta_kA)),\qquad
A~\hbox{ample on $X$,}~~\delta_k\le c{\log k\over k},~~c\ll 1,
$$
and to show that $Y_k=\pi_k(B_k)$ is a proper algebraic 
subvariety of~$X$ for $k$ large, under the assumption that $K_V$ is big. 
This does not look completely hopeless, 
since the statistics of curvature in the Morse inequalities do involve
currents for which the sets of poles depend only on the bigness of
$K_V$ and therefore project onto a proper subvariety $S$ of $X$ (see
the last step of the proof in section~$2$). It is not unreasonable to
think that a further analysis of the asymptotic behavior of sections,
e.g.\ through estimates of the Bergman kernel, might lead to such
results.

Even if the required property of the base locus cannot be obtained
directly, it would be enough, for a suitable irreducible analytic set
$Z\subset X_k^\GG$ contained in the base locus at some stage,
to construct non zero sections in
$$H^0(Z,\cO_{X_k^\GG}(m)_{|Z}\otimes\pi_k^*\cO(-m\delta_kA)_{|Z})$$
whenever $\pi_k(Z)=X$, and then to proceed inductively to cut-down
the base locus until one reaches some $Z'\subset Z$ with 
$\pi_k(Z')\subsetneq X$. Hence we have to estimate the cohomology groups
$H^0$ and $H^q$ not just on $X_k^\GG$, but also on all irreducible
subvarieties $Z\subset X_k^\GG$ such that $\pi_k(Z)=X$. We are not able
to do this in such a generality, but our method does provide
interesting results in this direction.

\claim (3.1) Theorem|Let $(X,V)$ be a compact directed $n$-dimensional
manifold, let \hbox{$r=\rank V$} and $F$ be a holomorphic line bundle on $X$. 
Fix an irreducible analytic set
$\smash{Z_{k_0}\subset X_{k_0}^\GG}$ or equivalently some
$\bC^*$-invariant set $\smash{Z'_{k_0}\subset J^{k_0}V}$, and assume that
$\smash{\pi_{k_0}(Z_{k_0})}=X$. For $k\gg k_0$,
denote by $Z_k\subset X_k^\GG$ the irreducible set corresponding to the
inverse image of $\smash{Z'_{k_0}}$ by the canonical morphism
\hbox{$J^kV\to J^{k_0}V$}. Let $h$ be an admissible metric on $V$, $h_F$ a
metric with analytic singularities on $F$ and
$$
\eqalign{
L_k&=\cO_{X_k^\GG}(1)\otimes
\pi_k^*\cO\Big({1\over kr}\Big(1+{1\over 2}+\ldots+{1\over k}\Big)F\Big),\cr
\eta&=\Theta_{K_V,\det h^*}+\Theta_{F,h_F},\qquad S=\Sing(\eta).\cr}
$$
Then for $m\gg k\gg k_0$ and $p_k=\dim Z_k=\dim Z_{k_0}+(k-k_0)r$ we have
$$\eqalign{
h^0(Z_k,&\cO(L_k^{\otimes m})_{|Z_k})\cr
&\ge {m^{p_k}\over p_k!}
{(\log k)^n\over n!}\deg_{X_k^\GG/X}(Z_k)\bigg(
\int_{X(\eta,\le 1)\ssm S}\eta^n-O((\log k)^{-1})\bigg)-o(m^{p_k})\cr}
$$
where $\deg_{X_k^\GG/X}(Z_k)=
\deg_{X_{k_0}^\GG/X}(Z_{k_0})\big({k_0!\over k!}\big)^r$
is the relative degree of $Z_k$ over $X$ with respect to the normalized
weighted ``K\"ahler metric'' $\omega_{a,r,p}$ introduced in~$(1.10)$.
\endclaim

We would also get similar upper and lower Morse bounds for the higher 
cohomology groups, provided that the sheaves $\smash{\cO_{X_k^\GG}(m)}$ are
twisted by the appropriate multiplier ideal sheaves $\cJ_{k,m}$ already 
described. The main trouble to proceed further in the analysis of the base
locus is that we have to take $k\gg k_0$ and that the $O(...)$ and $o(...)$
bounds depend on $Z_{k_0}$. Hence the newer sections can only be
constructed for higher and higher orders $k$, without any indication that
we can actually terminate the process somewhere, except possibly by
some extremely delicate uniform estimates which seem at present beyond 
reach.\medskip

\proof. The technique is a minor variation of what has been done in
section~2, hence we will only indicate the basic idea. Essentially
the $k$-jet of $f$ is no longer completely random, its projection onto
the first $k_0$ components $(\nabla^jf(0))_{1\le j\le k_0}$ is assigned 
to belong to some given analytic set $\smash{Z'_{k_0}\subset J^{k_0}}V$.
This means that in the curvature formula (2.15)
$$
g_k(z,x,u)={i\over 2\pi}\sum_{1\le s\le k}{1\over s}x_s
\sum_{i,j,\alpha,\beta}c_{ij\alpha\beta}(z)\,
u_{s\alpha}\overline u_{s\beta}\,dz_i\wedge d\overline z_j
$$
only the sum $\sum_{k_0<s\le k}$ is perfectly random. The partial sum
$\sum_{1\le s\le k_0}$ remains bounded, while the harmonic
series diverges as $\log k$. This implies that the ``non randomness'' of the
initial terms perturbs the estimates merely by bounded quantities, and
in the end, the expected value is still similar to (2.18), i.e.
$$
{\bf E}(g_k(z,\bu,\bu))={1\over kr}
\Big(1+{1\over 2}+\ldots+{1\over k}+O(1)\Big)\big(\Theta_{K_V,\det h^*}+
\Theta_{F,h_F}\big).
$$
Once we are there, the calculation of standard deviation and the other
estimates are just routine, and Theorem~3.1 follows again from Proposition
2.13 when we integrate the Morse integrals over~$Z_k$ instead of the
whole $k$-jet space $X_k^\GG$.\qed

Another possibility to analyze the base locus is to study the 
restriction maps
$$
\rho_{k,m}(x):H^0(X,E^\GG_{k,m}V^*\otimes\cO(-m\delta_kA))\to 
\big(E^\GG_{k,m}V^*\otimes\cO(-m\delta_kA)\big)_x\leqno(3.2)
$$
at generic points $x\in X$. If $\rho_{k,m}(x)$ can be shown to be surjective
at a generic point, then a fortiori the projection $Y_k=\pi_k(B_k)$ of the 
base locus does not contain $x$ and so $Y_k$ is a proper algebraic subvariety 
of~$X$. Now, proving the surjectivity of $\rho_{k,m}(x)$ could be done by
proving the vanishing of the $H^1$ group of our sheaf twisted by
the maximal ideal $\gm_{X,x}$. We cannot exactly reach such a precise
vanishing result, but Morse inequalities can be used to show that the 
$H^1$ groups do not grow too fast.

In fact assume that $A$ is an ample $\bQ$-divisor on $X$ which is chosen 
so small that $K_V\otimes\cO(-A)$ is still big. By our estimates, we 
can then take $\delta_k={1\over kr}(1+{1\over 2}+\ldots+{1\over k})$.
Pick a very ample divisor $G$ on $X$ and $n$ pencils
of sections $\sigma_{j,t}\in H^0(X,\cO(G))$, $1\le j\le n$, $t\in\bP^1_\bC$,
such that the divisors $\sigma_{j,t_j}(z)=0$ intersect transversally
at isolated points for generic choices of the parameters $t_j\in\bP^1_\bC$.
We select an admissible metric $h$ on $V$ which provides a strictly
positive curvature current on $K_V\otimes\cO(-A)$ and multiply it
by the additional weight factor $(e^\varphi)^{1/rm\delta_k}$ where 
$$\varphi(z)=\log\sum_{1\le j\le n}\prod_{t\in T_j}|\sigma_{j,t}(z)|^{2n}_{h_G}$$
and $T_j\subset \bP^1_\bC$ are generic finite subsets of given 
cardinality~$N$. The multiplier ideal sheaf of
$\varphi$ is precisely equal to the ideal $\cI_E$ of germs of
functions vanishing on a certain $0$-dimensional set 
$E=\{x_1,\ldots,x_s\}\subset X$ of cardinality $s=N^nG^n$. 
Also the resulting curvature form
$$
\eta=\Theta_{K_V,\det h^*}-\Theta_{A,h_A}+{1\over m\delta_k}dd^c\varphi
\ge \Theta_{K_V,\det h^*}-\Theta_{A,h_A}-{N\over m\delta_k}\Theta_{G,h_G}
$$
can be made to be strictly positive as a current provided that
$N\sim cm\delta_k$ with~$c\ll 1$. Then the corresponding multiplier
ideal sheaf of the induced hermitian metric on
$$
\cO_{X_k^\GG}(m)\otimes \pi_k^*\cO(-m\delta_k A)
$$
is the original multiplier sheaf $\cJ_{k,m}$ twisted 
by $\pi_k^*\cI_E$ above $x_j$, provided that the $x_j$ lie outside 
of $\Sing(V)$ and outside of the projection of the support $V(\cJ_{k,m})$. 
Consider the exact sequence
$$
\eqalign{
0&\lra \cO_{X_k^\GG}(m)\otimes \pi_k^*\cO(-m\delta_k A)\otimes\cJ_{k,m}
\otimes\pi_k^*\cI_E\cr
&\lra
\cO_{X_k^\GG}(m)\otimes \pi_k^*\cO(-m\delta_k A)\otimes\cJ_{k,m}\cr
&\lra \cO_{X_k^\GG}(m)\otimes \pi_k^*\cO(-m\delta_k A)\otimes\cJ_{k,m}
\otimes\pi_k^*(\cO_X/\cI_E) \lra 0.\cr}
$$
Its cohomology exact sequence yields an ``almost surjective arrow''
$$
H^0\big(\cO_{X_k^\GG}(m)\otimes
\pi_k^*\cO(-m\delta_k A)\otimes\cJ_{k,m}\big)
\lra\bigoplus_{1\le j\le s}\big(E^\GG_{k,m}V^*\otimes\cO(-m\delta_kA)\big)_{x_j},
$$
namely the image contains the kernel of the map
$$
\bigoplus_{1\le j\le s}\big(E^\GG_{k,m}V^*\otimes\cO(-m\delta_kA)\big)_{x_j}
\lra
H^1\big(\cO_{X_k^\GG}(m)\otimes
\pi_k^*\cO(-m\delta_k A)\otimes\cJ_{k,m}\otimes
\pi_k^*\cI_E\big).
$$
Now, we have a Morse upper bound
$$
h^1\big(\cO_{X_k^\GG}(m)\otimes
\pi_k^*\cO(-m\delta_k A)\otimes\cJ_{k,m}\otimes
\pi_k^*\cI_E\big)\le
{m^{n+kr-1}\over (n+kr-1)!}{(\log k)^n\over n!\,(k!)^r}O\big((\log k)^{-1}\big)
$$
since the $1$-index integral $\int_{X(\eta,1)}h^n$ is identically zero. At the
same time we have $s=N^nG^n\sim c'm^n(\log k)^n/k^n$, and 
it follows that 
$$\dim\bigoplus_{1\le j\le s}\big(E^\GG_{k,m}V^*\otimes
\cO(-m\delta_kA)\big)_{x_j}
\sim s{m^{kr-1}\over (kr-1)!(k!)^r}\sim
{c'm^{n+kr-1}\over (kr-1)!(k!)^r}{(\log k)^n\over k^n}.
$$
By selecting a suitable point $x_j$ and a trivial lower semi-continuity
argument we get the desired almost surjectivity.

\claim (3.3) Corollary|If $A$ is an ample $\bQ$-divisor on $X$ such that
$K_V\otimes\cO(-A)$ is big and $\delta_k=
{1\over kr}(1+{1\over 2}+\ldots+{1\over k})$, $r=\rank V$,
the restriction map
$$
\rho_{k,m}(x):H^0(X,E^\GG_{k,m}V^*\otimes\cO(-m\delta_kA))\to 
\big(E^\GG_{k,m}V^*\otimes\cO(-m\delta_kA)\big)_x
$$
has an image of dimension larger than $(1-O((\log k)^{-1}))\dim E^\GG_{k,m}V^*$
at a generic point $x\in X$ for $m\gg k\gg 1$.
\endclaim

Such a result puts an upper bound on the vanishing order that a generic
section may have on $X_k^\GG$ above a generic point of~$X$. Our hope is 
that one can then completely ``eliminate'' the base locus by taking 
vertical derivatives along the fibers of $J^kV\to X\,$; those derivations will
necessarily have some poles $\cO(pA)$ which we hope to get cancelled by 
the negative powers $\cO(-m\delta_kA)$. This strategy 
first devised by [Siu02, Siu04] has indeed been successful in some cases 
for the study of generic algebraic degeneracy (e.g.\ for hypersurfaces
of very large degree in $\bP^{n+1}_\bC$). This would work rather easily 
if the rough error term $O((\log k)^{-1})$ could be replaced e.g.\ by
$O(m^{-\varepsilon_k})$ in Corollary (3.3), but this is maybe too much 
to ask~for.

We finally discuss yet another approach. For this we have to introduce
invariant jet differentials along the lines of [Dem95]. In fact, to any 
directed manifold $(X,V)$ one can associate its tower of Semple 
$k$-jet spaces, which is a sequence of
directed pairs $(X_k,V_k)$ starting with $(X_0,V_0)=(X,V)$, together
with morhisms $\wt\pi_k:(X_k,V_k)\to (X_{k-1},V_{k-1})$. These spaces
are constructed inductively by putting $X_k=P(V_{k-1})$ and
$V_k=(\wt\pi_k)_*^{-1}(\cO_{X_k}(-1))$ where
$$
\cO_{X_k}(-1)\subset (\wt\pi_k)^*V_{k-1}\subset (\wt\pi_k)^*T_{X_{k-1}}
$$
is the tautological subbundle (cf.~[Dem95]). In the case where
$V$ is not a subbundle, we can first construct the
absolute tower $(\ol X_k,\ol V_k)$ by starting from $\ol V_0=T_X$, 
and then take $X_k$ to be the closure in $\ol X_k$ of the $k$-step $X'_k$ 
of the relative tower $(X'_k,V'_k)$ constructed 
over the dense Zariski open set $X'=X\ssm \Sing(V)$.
In this way, the tower $(X_k,V_k)$ is 
at least birationally well defined -- in such a birational context we 
can even assume that $X_k$ is
smooth after performing a suitable modification at each stage. Even if 
we start with $V=T_X$ (or an integrable subbundle $V\subset T_X$), the
$k$-jet lifting $V_k$ will not be integrable in general, the only
exception being when $\rank V_k=\rank V=1$. Now, if
$$
\pi_{k,0}=\wt\pi_k\circ\ldots\circ\wt\pi_1:X_k\to X_0=X,
$$
it is shown in [Dem95] that  the direct image sheaf
$$
\pi_{k,0}\cO_{X_k}(m):=E_{k,m}V^*\subset E_{k,m}^\GG V^*
$$
consists of algebraic differential operators $P(f^{(j)}_{j\le k})$ which 
satisfy the invariance property
$$
P((f\circ\varphi)^{(j)}_{j\le k})=
(\varphi')^mP(f^{(j)}_{j\le k})\circ \varphi
$$
when $\varphi\in \bG_k$ is in the group of $k$-jets of biholomorphisms
$\varphi:(\bC,0)\to(\bC,0)$. Since we already assume $\bC^*$ invariance,
it is enough to require invariance by the nilpotent subgroup
$\bG'_k\subset\bG_k$ of $k$-jets tangent to identity. The group
$\bG'_k$ is a semi-direct product of additive groups $(\bC,+)$ consisting
of biholomorphisms $\tau_{j,a}:t\mapsto t+at^j+O(t^{j+1})$, $2\le j\le k$,
$a\in\bC$. In this tower, the biholomorphisms $\tau_{k,a}$ actually
generate a normal subgroup of $\bG'_k$, and we have
$\bG'_k/\{\tau_{k,a}\}\simeq\bG'_{k-1}$. Now, assume that we have found
a section
$$
P\in H^0(X,E^\GG_{k,m}V^*\otimes\cO(-m\delta_kA))
$$
for some ample $\bQ$-divisor $A$ on $A$. Then we have an expansion
$$
P_a(f^{(j)}_{j\le k}):=P((f\circ\tau_{k,a})^{(j)}_{j\le k})=
\sum_{0\le s\le m/k}a^sP_s(f^{(j)}_{j\le k})
$$
and the highest  non zero term $P_s$ is $\{\tau_{k,a}\}$-invariant of 
weighted degree $m-(k-1)s\,$; this comes from the fact that the homothety 
$h_\lambda(t)=\lambda t$ satisfies 
$$\tau_{k,a}\circ h_\lambda=h_\lambda\circ\tau_{k,a\lambda^{k-1}}.$$
Then it makes sense to look at the action of 
$\{\tau_{k-1,a}\}$ on $P_s$,
and proceeding inductively we reach a non zero $\bG'_k$-invariant
(and thus $\bG_k$-invariant) polynomial
$$
Q\in H^0(X,E_{k,m'}V^*\otimes\cO(-m\delta_kA))
$$
of degree $m'\le m$ (and possibly of order $k'\le k$ but we can still
consider it to be of order $k$). By raising $Q$ to some power $p$ and using
the $\bQ$-ampleness of $A$, we obtain a genuine integral section
$$
Q^p\sigma_A^{p(m-m')\delta_k}\in H^0(X,E_{k,pm'}V^*\otimes\cO(-pm'\delta_kA)).
$$

\claim (3.4) Corollary|Let $(X,V)$ be a projective directed manifold such that
$K_V$ is big, and $A$ an ample $\bQ$-divisor on $X$ such that 
$K_V\otimes\cO(-A)$ is still big. Then, if we put $\delta_k={1\over kr}
(1+{1\over 2}+\ldots+{1\over k})$, $r=\rank V$, the space of global 
invariant jet differentials
$$
H^0(X,E_{k,m}V^*\otimes\cO(-m\delta_kA))
$$
has $($many$)$ non zero sections for $m\gg k\gg 1$.
\endclaim

If we have a directed projective variety $(X,V)$ with $K_V$ big, we conclude
that there exists $k\ge 1$ and a proper analytic set $Z\subset X_k$ such that
all entire curves have the image of their $k$-jet $f_{[k]}(\bC)$ contained
in $Z$. Let $Z'$ be an irreducible component of $Z$ such that 
$\pi_{k,0}(Z')=X$ (if $\pi_{k,0}(Z')\subsetneq X$ there is nothing more
to do). Consider the linear subspace
$V'=\overline{T_{Z'\ssm Z''}\cap V_k}$ where $Z''\subset Z'$ is chosen 
such that $Z'\ssm Z''$ is non singular and the intersection 
$T_{Z'\ssm Z''}\cap V_k$ is a subbundle of $T_{Z'\ssm Z''}$. 
If $f_{[k]}(\bC)$ is not
contained identically in $Z''$, then the curve $g=f_{[k]}$ is tangent
to $(Z',V')$. On the other hand, if $f_{[k]}(\bC)\subset Z''$ we can replace 
$Z'$ by $Z''$ and argue inductively on $\dim Z'$. What we have gained 
here is that we have
replaced the initial directed space $(X,V)$ with another one $(Z',V')$
such that $\rank V'<\rank V$, and we can try to argue by induction on
$r=\rank V$.

Observe that the generalized Green-Griffiths conjecture
is indeed trivial for $r=1$ (assuming $K_V=\cO(V^*)$ big)$\,$: in fact we get 
in this case a non zero section $P\in H^0(X,V^{*\otimes k}
\otimes\cO(-A))$ for some $k\gg 1$ and so $P(f)\cdot (f')^k$ must vanish 
for every entire curve $f:(\bC,T_\bC)\to (X,V)$. Therefore $f(\bC)\subset
Y:=\{P(z)=0\}\subsetneq X$. The main difficulty in this inductive
approach is that when we start with $(X,V)$ with $K_V$ big, it seems to 
be very hard to say anything about $K_{V'}$ on $(Z',V')$. Especially,
the singularities of $Z'$ and $V'$ do not seem to be under control. The only
hope would be to have enough control on the sections cutting out $Z'$, and this
requires anyway to understand much more precisely the behavior and 
vanishing order of generic sections
$P\in H^0(X,E_{k,m}V^*\otimes\cO(-m\delta_kA))$. One could try in this context
to take $A$ to approach the positive part in the Zariski decomposition
of~$K_V$, in such a way that the sections $P$ do not have much space to
move around statistically.

\section{References}
\medskip

{\eightpoint

\bibitem[Ber10]&B\'erczi, G.:& Thom polynomials and the
  Green-Griffiths conjecture.& Manuscript Math.\ Institute Oxford, May 2010&

\bibitem[Blo26]&Bloch, A.:& Sur les syst\`emes de fonctions uniformes
  satisfaisant \`a l'\'equation d'une vari\'et\'e alg\'ebrique dont
  l'irr\'egularit\'e d\'epasse la dimension.& J.\ de Math., {\bf 5}
  (1926), 19--66&

\bibitem[Bog79]&Bogomolov, F.A.:& Holomorphic tensors and vector
  bundles on projective varieties.& Math.\ USSR Izvestija {\bf 13}/3
  (1979), 499--555&

\bibitem[Bon93]&Bonavero, L.:& In\'egalit\'es de morse holomorphes 
  singuli\`eres.&
  Acad.\ Sci.\ Paris S\'er.~I Math.\ {\bf 317} (1993), 1163–-1166, and:  
  J.\ Geom.\ Anal.\ {\bf 8} (1998), 409-–425&

\bibitem[Bro78]&Brody, R.:& Compact manifolds and hyperbolicity.&
  Trans.\ Amer.\ Math.\ Soc.\ {\bf 235} (1978), 213--219&

\bibitem[BG77]&Brody, R., Green, M.:& A family of smooth hyperbolic
  surfaces in $\bP^3$.& Duke Math.\ J.\ {\bf 44} (1977), 873--874&

\bibitem[Cle86]&Clemens, H.:& Curves on generic hypersurfaces.&
  Ann.\ Sci.\ \'Ec. Norm.\ Sup.\ {\bf 19} (1986) 629--636, Erratum:
  Ann.\ Sci.\ \'Ec. Norm.\ Sup.\ {\bf 20} (1987) 281&

\bibitem[CG76]&Cowen, M., Griffiths, P.:& Holomorphic curves and
  metrics of negative curvature.& J.~Analyse Math.\ {\bf 29} (1976),
  93--153&

\bibitem[Dem85]&Demailly, J.-P.:& Champs magn\'etiques et in\'egalit\'es 
  de Morse pour la $d''$-cohomologie.& Ann.\ Inst.\ Fourier (Grenoble), 
  {\bf 35} (1985), 189--229&

\bibitem[Dem95]&Demailly, J.-P.:& Algebraic criteria for Kobayashi
  hyperbolic projective varieties and jet differentials.& AMS Summer
  School on Algebraic Geometry, Santa Cruz 1995, Proc.\ Symposia in
  Pure Math., ed.\ by J.~Koll\'ar and R.~Lazarsfeld, 76p&

\bibitem[Dem97]&Demailly, J.-P.:& Vari\'et\'es hyperboliques et
  \'equations diff\'erentielles alg\'ebriques.& Gaz.\ Math.\ {\bf 73}
  (juillet 1997), 3--23&

\bibitem[Dem10a]&Demailly, J.-P.:& Holomorphic Morse inequalities and
  asymptotic cohomology groups: a tribute to Bernhard Riemann.&
  arXiv: math.CV/1003.5067&

\bibitem[Dem10b]&Demailly, J.-P.:& A converse to the
  Andreotti-Grauert theorem.& Manuscript Inst.\ Fourier Grenoble, 
  October 2010&

\bibitem[DEG00]&Demailly, J.-P., El Goul, J.:& Hyperbolicity of
  generic surfaces of high degree in projective 3-space.&
  Amer.\ J.\ Math.\ {\bf 122} (2000), 515-\-546&

\bibitem[Div09]&Diverio, S.:& Existence of global invariant jet
  differentials on projective hypersurfaces of high degree.& Math.\
  Ann.\ {\bf 344} (2009), 293--315&

\bibitem[DMR10]&Diverio, S., Merker, J., Rousseau, E.:& Effective
  algebraic degeneracy.& Invent.\ Math.\ {\bf 180} (2010), 161--223&

\bibitem[DT10]&Diverio, S., Trapani, S.:& A remark on the codimension
  of the Green-Griffiths locus of generic projective hypersurfaces of
  high degree.&To appear in J.\ Reine Angew.\ Math&

\bibitem [Fuj94]&Fujita, T.:& Approximating Zariski decomposition of
  big line bundles.& Kodai Math.\ J.\ {\bf 17} (1994), 1--3&

\bibitem[GG79]&Green, M., Griffiths, P.:& Two applications of
  algebraic geometry to entire holomorphic mappings.& The Chern
  Symposium 1979, Proc.\ Internal.\ Sympos.\ Berkeley, CA, 1979,
  Springer-Verlag, New York (1980), 41--74&

\bibitem[Kob70]&Kobayashi, S.:& Hyperbolic manifolds and holomorphic
  mappings.& Marcel Dekker, New York (1970)&

\bibitem[Kob76]&Kobayashi, S.:& Intrinsic distances, measures and
  geometric function theory.& Bull.\ Amer.\ Math.\ Soc.\ {\bf 82}
  (1976), 357--416&

\bibitem[KobO75]&Kobayashi, S., Ochiai, T.:& Meromorphic mappings into
  compact complex spaces of general type.& Invent.\ Math.\ {\bf 31}
  (1975), 7--16&

\bibitem[Lang86]&Lang, S.:& Hyperbolic and Diophantine analysis.&
  Bull.\ Amer.\ Math.\ Soc.\ {\bf 14} (1986) 159--205&

\bibitem[Lang87]&Lang, S.:& Introduction to complex hyperbolic
  spaces.& Springer-Verlag, New York (1987)&

\bibitem[McQ96]&McQuillan, M.:& A new proof of the Bloch conjecture.&
  J.\ Alg.\ Geom.\ {\bf 5} (1996), 107--117&

\bibitem[McQ98]&McQuillan, M.:& Diophantine approximation and foliations.& 
  Inst.\ Hautes \'Etudes Sci.\ Publ.\ Math.\ {\bf 87} (1998), 121--174&

\bibitem[McQ99]&McQuillan, M.:& Holomorphic curves on hyperplane sections 
  of $3$-folds.& Geom.\ Funct.\ Anal.\ {\bf 9} (1999), 370--392&

\bibitem[Mer08]&Merker, J.:& Jets de Demailly-Semple d’ordres 4 et 5 
  en dimension 2.& arXiv: math.AG/0710.2393, 
  Int.\ J.\ Contemp.\ Math.\ Sciences {\bf 3} (2008) 861--933&

\bibitem[Mer09]&Merker, J.:& Low pole order frames on vertical jets of
  the universal hypersurface.& arXiv: math.AG/0805.3987, Ann.\ Inst.\
  Fourier (Grenoble) {\bf 59} (2009) 1077–-1104&

\bibitem[Mer10]&Merker, J.:& Complex projective hypersurfaces of
  general type: toward a conjecture of Green and Griffiths.&
  Manuscript \'Ec.\ Norm.\ Sup.\ Paris, May 2010, arXiv: math.AG/1005.0405&

\bibitem[Pau08]&P\u{a}un, M.:& Vector fields on the total space of 
  hypersurfaces in the projective space and hyperbolicity.& Math.\ Ann.\ 
  {\bf 340} (2008) 875--892&

\bibitem[Rou06a]&Rousseau, E.:& \'Etude des jets de Demailly-Semple en 
  dimension 3.& Ann.\ Inst.\ Fourier (Grenoble) {\bf 56} (2006) 397--421&

\bibitem[Rou06b]&Rousseau, E.:& \'Equations diff\'erentielles sur les 
  hypersurfaces de $\bP^4$.& J.\ Math.\ Pures Appl.\ {\bf 86} (2006)
  322--341&

\bibitem[Rou07]&Rousseau, E.:& Weak analytic hyperbolicity of generic 
  hypersurfaces of high degree in $\bP^4$.& Annales Fac.\ Sci.\ Toulouse
  {\bf 16} (2007), 369--383&

\bibitem[Siu97]&Siu, Y.T.:& A proof of the general schwarz lemma using
  the logarithmic derivative lemma.& Communication personnelle, avril
  1997&

\bibitem[Siu02]&Siu, Y.T.:& Some recent transcendental techniques in
  algebraic and complex geometry.& Proceedings of the International
  Congress of Mathematicians, Vol.~I (Beijing, 2002), Higher
  Ed.\ Press, Beijing, 2002, 439--448&

\bibitem[Siu04]&Siu, Y.T.:& Hyperbolicity in complex geometry.& The 
  legacy of Niels Henrik Abel, Springer, Berlin, 2004, 543--566&

\bibitem[SY96a]&Siu, Y.T., Yeung, S.K.:& Hyperbolicity of the
  complement of a generic smooth curve of high degree in the complex
  projective plane.& Invent.\ Math.\ {\bf 124} (1996), 573--618&

\bibitem[SY96b]&Siu, Y.T., Yeung, S.K.:& Defects for ample divisors
  of Abelian varieties, Schwarz lemma and hyperbolic surfaces of low
  degree.& Preprint (pr\'epublication, automne 1996)&

\bibitem[Tra95]&Trapani, S.:& Numerical criteria for the positivity
  of the difference of ample divisors.& Math.\ Z.\ {\bf 219} (1995),
  387--401&

\bibitem[Voi96]&Voisin, C.:& On a conjecture of Clemens on rational
  curves on hypersurfaces.& J.\ Diff.\ Geom.\ {\bf 44} (1996),
  200--213. Correction: J.\ Diff.\ Geom.\ {\bf 49} (1998), 601-\-611&

}
\vskip20pt
\noindent
(version of November 26, 2010, printed on \today)

\end